\documentclass{article}

\usepackage{amsmath}
\usepackage{amsthm}
\usepackage{amssymb}
\usepackage{amsfonts}
\usepackage{mathrsfs}

\usepackage{subeqnarray}
\usepackage{bm}
\usepackage{tikz}
\usepackage{dsfont}
\usepackage{upgreek}
\usepackage{booktabs}
\usepackage{fancyhdr}
\usepackage{multirow}
\usepackage{multicol}
\usepackage{float}
\usepackage[inline]{enumitem}
\usepackage{cases}
\usepackage{xurl}

\usepackage[ruled, vlined, linesnumbered, inoutnumbered, algo2e]{algorithm2e}

\usepackage{bbding}
\usepackage{caption}

\usepackage[mathlines]{lineno}

\usepackage[top = 2.5cm, bottom = 2.5cm, left = 3cm, right = 3cm]{geometry}

\usepackage{graphics, graphicx, epsf, epstopdf, subfig}
\graphicspath{{./Figures/}}
\ifpdf
\DeclareGraphicsExtensions{.eps,.pdf,.png,.jpg}
\else
\DeclareGraphicsExtensions{.eps}
\fi

\usepackage[breaklinks, colorlinks, linkcolor = black, citecolor = black, urlcolor = black]{hyperref}

\usepackage{xcolor}
\usepackage{array}
\usepackage{placeins}
\usepackage{pgfplots}
\pgfplotsset{compat = 1.18}
\usetikzlibrary{calc, arrows.meta, positioning}

\providecommand{\U}[1]{\protect\rule{.1in}{.1in}}

\newtheorem{theorem}{Theorem}[section]
\newtheorem{corollary}[theorem]{Corollary}
\newtheorem{lemma}[theorem]{Lemma}
\newtheorem{proposition}[theorem]{Proposition}
\newtheorem{definition}[theorem]{Definition}

\newtheorem{assumption}{Assumption}

\newtheorem{example}{Example}

\numberwithin{equation}{section}

\newcommand{\diag}{\mathrm{diag}}
\newcommand{\Diag}{\mathrm{Diag}}
\newcommand{\conv}{\mathrm{conv}}
\newcommand{\tr}{\mathrm{tr}}
\newcommand{\proj}{\mathrm{Proj}}

\newcommand{\grad}{\mathrm{grad}}

\newcommand{\bN}{\mathbb{N}}
\newcommand{\bR}{\mathbb{R}}

\newcommand{\cB}{\mathcal{B}}
\newcommand{\cC}{\mathcal{C}}
\newcommand{\cE}{\mathcal{E}}
\newcommand{\cI}{\mathcal{I}}
\newcommand{\cJ}{\mathcal{J}}

\newcommand{\cS}{\mathcal{S}}
\newcommand{\cT}{\mathcal{T}}
\newcommand{\cU}{\mathcal{U}}

\newcommand{\sB}{\mathscr{B}}
\newcommand{\sH}{\mathscr{H}}

\newcommand{\sff}{\mathsf{f}}

\newcommand{\sfm}{\mathsf{m}}
\newcommand{\sfr}{\mathsf{r}}

\newcommand{\sfE}{\mathsf{E}}
\newcommand{\sfF}{\mathsf{F}}
\newcommand{\sfG}{\mathsf{G}}

\newcommand{\sfJ}{\mathsf{J}}
\newcommand{\sfL}{\mathsf{L}}
\newcommand{\sfM}{\mathsf{M}}
\newcommand{\sfN}{\mathsf{N}}
\newcommand{\sfQ}{\mathsf{Q}}
\newcommand{\sfR}{\mathsf{R}}
\newcommand{\sfS}{\mathsf{S}}
\newcommand{\sfT}{\mathsf{T}}
\newcommand{\sfU}{\mathsf{U}}
\newcommand{\sfV}{\mathsf{V}}
\newcommand{\sfW}{\mathsf{W}}
\newcommand{\sfY}{\mathsf{Y}}

\newcommand{\dE}{\mathds{E}}
\newcommand{\dJ}{\mathds{J}}
\newcommand{\dN}{\mathds{N}}
\newcommand{\dP}{\mathds{P}}
\newcommand{\dQ}{\mathds{Q}}
\newcommand{\dW}{\mathds{W}}

\newcommand{\rmD}{\mathrm{D}}

\newcommand{\Rr}{\mathbb{R}^{r}}
\newcommand{\Rd}{\mathbb{R}^{d}}

\newcommand{\Rdr}{\mathbb{R}^{d \times r}}
\newcommand{\Rs}{\mathbb{R}^{s}}
\newcommand{\Rss}{\mathbb{R}^{s \times s}}

\newcommand{\Rdn}{\mathbb{R}^{d \times n}}

\newcommand{\Odr}{\mathcal{O}^{d, r}}
\newcommand{\Ods}{\mathcal{O}^{d, s}}

\newcommand{\bfone}{\mathbf{1}}

\newcommand{\zz}{^{\top}}
\newcommand{\inv}{^{-1}}
\newcommand{\ff}{_{\mathrm{F}}}
\newcommand{\fs}{^2_{\mathrm{F}}}

\newcommand{\ucirc}{^{\circ}}
\newcommand{\uast}{^{\ast}}
\newcommand{\ustar}{^{\star}}
\newcommand{\udiamond}{^{\diamond}}

\newcommand{\lbest}{_{\mathrm{best}}}
\newcommand{\lrec}{_{\mathrm{rec}}}
\newcommand{\ltest}{_{\mathrm{test}}}
\newcommand{\utest}{^{\mathrm{test}}}
\newcommand{\ltrain}{_{\mathrm{train}}}
\newcommand{\utrain}{^{\mathrm{train}}}
\newcommand{\lbad}{_{\mathrm{bad}}}
\newcommand{\ubad}{^{\mathrm{bad}}}

\newcommand{\lres}{_{\mathrm{res}}}

\newcommand{\lpca}{_{\mathrm{pca}}}

\newcommand{\lrmR}{_{\mathrm{R}}}

\newcommand{\lDSGM}{_{\mathrm{DSGM}}}
\newcommand{\lRSG}{_{\mathrm{RSG}}}

\newcommand{\MSE}{\mathrm{MSE}}
\newcommand{\AMR}{\mathrm{AMR}}
\newcommand{\AV}{\mathrm{AV}}
\newcommand{\MD}{\mathrm{MD}}
\newcommand{\SR}{\mathrm{SR}}

\newcommand{\esssup}{\mathop{\mathrm{ess\, sup}}}

\newcommand{\dkh}[1]{\left(#1\right)}
\newcommand{\hkh}[1]{\left\{#1\right\}}
\newcommand{\fkh}[1]{\left[#1\right]}
\newcommand{\jkh}[1]{\left\langle#1\right\rangle}
\newcommand{\norm}[1]{\left\|#1\right\|}
\newcommand{\abs}[1]{\left\lvert #1\right\rvert}
\newcommand{\floor}[1]{\left\lfloor #1\right\rfloor}
\newcommand{\ceil}[1]{\left\lceil #1\right\rceil}

\newcommand{\lport}{_{\mathrm{port}}}

\allowdisplaybreaks
\DeclareMathOperator*{\argmin}{arg\, min}
\DeclareMathOperator*{\argmax}{arg\, max}

\title{Infinity-Wasserstein Distributionally Robust Principal Component Analysis: Exact Reformulation, Finite-Sample Guarantee, and Efficient Algorithm\thanks{This work is supported by National Key R\&D Program of China (2023YFA1009300), National Natural Science Foundation of China (12125108, 12631014, 12288201), and RGC grant JLFS/P-501/24 for the CAS AMSS-PolyU Joint Laboratory in Applied Mathematics, CAS-Croucher Funding Scheme for Joint Laboratories.}}

\author{
	Zhouxing Luo\thanks{Institute of Computational Mathematics and Scientific/Engineering Computing, Academy of Mathematics and Systems Science, Chinese Academy of Sciences, Beijing, China (\href{mailto:luozhouxing18@mails.ucas.ac.cn}{luozhouxing18@mails.ucas.ac.cn}).}
	\and
	Lei Wang\thanks{Department of Applied Mathematics, The Hong Kong Polytechnic University, Hung Hom, Kowloon, Hong Kong, China (\href{mailto:wlkings@lsec.cc.ac.cn}{wlkings@lsec.cc.ac.cn}).}
	\and
	Xin Liu\thanks{Institute of Computational Mathematics and Scientific/Engineering Computing, Academy of Mathematics and Systems Science, Chinese Academy of Sciences, Beijing, China (\href{mailto:liuxin@lsec.cc.ac.cn}{liuxin@lsec.cc.ac.cn}).}
}

\date{}

\begin{document}


\maketitle

\begin{abstract}
	Principal component analysis (PCA) is a fundamental and ubiquitous technique for dimension reduction and representation learning, yet its empirical formulation remains vulnerable to distributional uncertainty. 
	While distributionally robust optimization (DRO) offers a natural remedy, existing order-2 Wasserstein models of PCA either preserve the nominal PCA subspace or lack an exact reformulation. 
	This paper demonstrates how the infinity-Wasserstein distance can reconcile these two desirable properties. 
	For the resulting DRO model of PCA, we derive an exact finite-dimensional reformulation and explicitly identify a worst-case distribution attaining the inner maximum. 
	Crucially, the reformulated objective retains samplewise residual information beyond the empirical covariance matrix, revealing a mechanism through which distributional robustness can reshape the learned principal subspace. 
	Under the mild assumption that the true distribution has a finite fourth moment, we further establish a uniform finite-sample performance guarantee for the proposed model. 
	An iteratively reweighted majorization–minimization algorithm is also developed to solve the resulting nonconvex and nonsmooth optimization problem, with a rigorous analysis of its convergence properties. 
	Numerical experiments substantiate the computational efficiency of the proposed algorithm and the practical effectiveness of the infinity-Wasserstein model in downstream applications. 
\end{abstract}


\section{Introduction}

Principal component analysis (PCA) \cite{Jolliffe2002Principal} is among the most widely used tools for dimension reduction in operations research, machine learning, and data science. 
By projecting high-dimensional observations onto a low-dimensional linear subspace, PCA provides a parsimonious representation of complex data and is routinely embedded in broader pipelines for forecasting, risk analysis, statistical learning, and optimization \cite{Fan2008High,Fan2013Large,Jiang2026Optimized,Kao2014Directed}. 
In most applications, however, the true distribution governing future observations is unknown and must be inferred from a finite sample. 
Classical PCA replaces the true distribution with its empirical counterpart, yielding an optimal solution given by the principal eigenspace of the empirical covariance matrix. 
Although this plug-in principle is computationally attractive, sampling error, data contamination, and distributional perturbations may distort the estimated covariance structure, which creates a substantial gap between in-sample fit and out-of-sample performance. 
The resulting discrepancy can, in turn, propagate to downstream prediction and optimization decisions that are built on the estimated subspace.

Distributionally robust optimization (DRO) offers a natural framework to address this uncertainty by optimizing against a family of distributions surrounding the empirical distribution \cite{Blanchet2019Quantifying,Delage2010Distributionally,Esfahani2018Data}. 
For PCA, the ambiguity set captured by Wasserstein distances is especially appealing because it describes distributional perturbations directly in the data space and therefore preserves the geometric interpretation of reconstruction error. 
The associated modeling challenge is more subtle than it first appears. 
A useful DRO formulation must change the selected principal subspace when the empirical distribution is unreliable, but it should also admit a sufficiently explicit characterization to support practical computation and finite-sample analysis. 
Existing Wasserstein formulations for PCA expose a tension between these two requirements.

Vu et al. \cite{Vu2022Distributionally} propose a distributionally robust model for fair PCA, in which the ambiguity set is quantified by the Gelbrich distance. 
For Gaussian distributions, the Gelbrich distance coincides with the order-$2$ Wasserstein distance. 
Wang et al. \cite{Wang2025Enhancing} further study distributionally robust PCA under the standard order-$2$ Wasserstein distance and derive an exact reformulation of the resulting min-max problem. 
For classical PCA without additional regularization, however, their characterization reduces the robust objective to a monotone transformation of the empirical reconstruction error. 
The nominal problem and its distributionally robust counterpart share exactly the same set of optimal solutions. 
This limitation identifies a basic modeling difficulty that the isotropic order-$2$ Wasserstein ambiguity does not by itself produce a new robust principal subspace in classical PCA.

To address this issue, Xu et al.~\cite{Xu2026Distributionally} replace the isotropic transportation cost with a weighted order-$2$ Wasserstein geometry. 
The anisotropic transport metric alters the relative ranking of candidate subspaces and generates an optimal subspace that differs from its nominal counterpart. 
However, this additional modeling flexibility comes at the expense of tractability, as the corresponding distributionally robust PCA problem does not, in general, admit an exact finite-dimensional reformulation. 
Xu et al.~\cite{Xu2026Distributionally} therefore resort to an upper-bound surrogate that combines the square root of the empirical reconstruction error with a spectral exposure term involving the largest eigenvalue of a matrix variable. 
To solve this surrogate problem, they further develop a heuristic procedure by searching over a structured finite family of candidate subspaces, for which no theoretical analysis is provided.

These observations leave open the fundamental question of whether one can construct a DRO model for PCA that induces a genuinely nontrivial robustification effect while still admitting an exact finite-dimensional reformulation. 
Ideally, such a model should also enjoy meaningful finite-sample performance guarantees and be amenable to efficient computation with rigorous convergence properties. 
In the present work, we answer this question affirmatively by devising a distributionally robust PCA model based on the infinity-Wasserstein distance that fulfills all of these requirements. 
The remainder of this paper is devoted to establishing its structural properties, statistical guarantees, and algorithmic foundations.

\subsection{Related Literature}

Recently, DRO has emerged as a powerful modeling paradigm for optimization under uncertainty, particularly in settings where the probability distributions of random parameters are only partially known. 
Interested readers are referred to \cite{Kuhn2025distributionally,Rahimian2022Frameworks} for a comprehensive overview of DRO. 
Beyond the aforementioned studies on DRO models of PCA \cite{Vu2022Distributionally,Wang2025Enhancing,Xu2026Distributionally}, we briefly review several additional strands of research that are relevant to this work.

Wasserstein DRO leverages the geometry of optimal transport to form ambiguity sets around an observed distribution \cite{Blanchet2019Quantifying,Esfahani2018Data,ShafieezadehAbadeh2019Regularization}. 
The order of Wasserstein distances plays a fundamental role in determining how transportation costs are aggregated across samples. 
In contrast to finite-order Wasserstein distances that constrain an average transportation cost, the infinity-Wasserstein distance imposes a uniform bound on the displacement of paired observations under an admissible joint distribution. 
This structural property explains why infinity-Wasserstein DRO is closely connected to sample robust optimization and adversarial perturbation models, while at the same time retaining a distributional interpretation. 
Early work by Xie~\cite{Xie2020Tractable} derives tractable convex reformulations for two-stage distributionally robust linear programs over infinity-Wasserstein ambiguity sets. 
From a statistical perspective, Nguyen et al.~\cite{Nguyen2020Distributionally} investigate distributionally robust nonparametric conditional estimation based on infinity-Wasserstein distances with finite-sample guarantees. 
A series of systematic studies by Bertsimas et al. \cite{Bertsimas2018Robust,Bertsimas2022Technical,Bertsimas2023Data} shows that sample robust optimization admits a natural interpretation in terms of infinity-Wasserstein ambiguity sets. 
Building on the same structural insight, Xie~\cite{Xie2021Distributionally} develops infinity-Wasserstein models for distributionally robust chance-constrained optimization. 
To alleviate the computational challenges incurred by infinity-Wasserstein distances, Wang et al.~\cite{Wang2024Regularization} integrate phi-divergence regularization into the underlying formulation. 
Recent results \cite{Gao2024Wasserstein,Zhang2025Short} further situate infinity-Wasserstein DRO within a unified framework by extending strong-duality and variation-regularization principles beyond finite-order transportation costs.

A complementary line of literature focuses on the robust estimation of covariance and precision matrices through DRO models. 
At a broader level, moment-based ambiguity sets, such as those introduced by Delage and Ye~\cite{Delage2010Distributionally} and Wiesemann et al.~\cite{Wiesemann2014Distributionally}, provide protections against uncertainty in moment information. 
Within the Wasserstein framework, Nguyen et al.~\cite{Nguyen2022Distributionally} develop a distributionally robust counterpart of Gaussian maximum-likelihood estimation. 
Blanchet and Si~\cite{Blanchet2019Optimal} further analyze the uncertainty size for this estimator, revealing that loss-optimal asymptotic calibration may differ from the canonical square-root scaling associated with distributional coverage. 
More recently, Yue et al.~\cite{Yue2026Geometric} formulate distributionally robust covariance estimation directly over the space of covariance matrices and elucidate how the geometry of ambiguity sets yields spectral shrinkage. 
These studies demonstrate that DRO can serve as a regularization mechanism for covariance estimation. 
By contrast, our approach incorporates distributional robustness directly into principal subspace identification, thereby preserving information that cannot be fully captured by the empirical covariance matrix alone.

\subsection{Contributions}

In this paper, we conduct a systematic exploration of a distributionally robust PCA model based on the infinity-Wasserstein distance. 
Our main contributions are summarized as follows.

An exact reformulation of the proposed model is first derived, which reduces the original infinite-dimensional min-max problem to a finite-dimensional minimization problem. 
The resulting representation identifies samplewise residual norms as a key source of distributional robustness, retaining information that cannot be captured by the covariance matrix alone. 
In addition, we explicitly construct a worst-case distribution attaining the inner maximum of the DRO model. 
These results lay the structural foundation for the subsequent analysis and algorithm design.

We further establish the finite-sample performance guarantee for the proposed model. 
Specifically, the DRO objective function is shown to provide an upper bound on the true objective with high probability, uniformly over all feasible points. 
Our analysis only relies on a finite fourth moment of the true distribution, thereby accommodating certain heavy-tailed distributions. 
This uniform control, in turn, ensures that the out-of-sample certificate continues to hold for data-depend choices.

The exact reformulation of the infinity-Wasserstein model leads to a nonconvex and nonsmooth optimization problem over the Stiefel manifold. 
By exploiting the underlying problem structure, we develop an iteratively reweighted majorization-minimization algorithm, whose iterations involve only the computation of a weighted sample mean and the solution of a weighted PCA problem. 
Despite its computational simplicity, the proposed algorithm enjoys rigorous theoretical guarantees, including the Clarke stationarity of accumulation points and an iteration-complexity bound.

Numerical experiments demonstrate the significant computational advantage of our algorithm over the general-purpose Riemannian optimization methods considered. 
Beyond algorithmic performance, the practical effectiveness of our DRO model is investigated in two real-world applications, namely, image classification and portfolio selection. 
The numerical findings indicate that more reliable principal subspace estimation not only enhances the quality of learned representations but also contributes to improved out-of-sample performance in downstream tasks.

\subsection{Organization}

The rest of this paper proceeds as follows. 
Section~\ref{sec:model_analysis} introduces the infinity-Wasserstein model of PCA and investigates its structural properties. 
We establish in Section~\ref{sec:statistical_properties}  the finite-sample performance guarantee of the proposed model. 
Section~\ref{sec:algorithm} is dedicated to designing an efficient algorithm to solve the exact reformulation, followed by convergence analysis in Section~\ref{sec:convergence}. 
Numerical results are presented in Section~\ref{sec:numerical} to examine the performance of our algorithm and model. 
Finally, this paper concludes with key insights and directions for future research in Section~\ref{sec:conclusion}.

\section{The Infinity-Wasserstein Model}

\label{sec:model_analysis}

This section aims to develop an infinity-Wasserstein DRO model for PCA. 
We derive an exact finite-dimensional reformulation of the resulting infinite-dimensional min-max problem. 
This transformation introduces a center variable whose structural properties are characterized subsequently. 
Building on the exact reformulation, we construct an explicit worst-case distribution and further clarify the connections and distinctions between the proposed formulation and existing order-$2$ Wasserstein DRO models. 

\subsection{Basic Notation}

Before presenting our model in detail, we first introduce some basic notation. 
The sets of real and natural numbers are denoted by $\bR$ and $\bN$, respectively. 
We use $\jkh{\cdot, \cdot}$ to represent the Euclidean inner product for vectors and its Frobenius counterpart for matrices. 
The notations $\norm{\cdot}_2$ and $\norm{\cdot}\ff$ stand for the Euclidean and Frobenius norms, respectively. 
We denote by $I_r$ the identity matrix of dimension $r$. 
The operator $\tr (\cdot)$ returns the trace of a given matrix. 
We adopt $\diag (\cdot)$ and $\Diag (\cdot)$ as the operations of extracting the diagonal entries of a given matrix as a vector and constructing a diagonal matrix from a given vector, respectively.

\subsection{Model Construction}

Let $\xi \in \Rd$ be a $d$-dimensional random vector governed by a probability distribution $\dP\uast$. 
To reduce the dimension of $\xi$ from $d$ to $r$, PCA attempts to seek an $r$-dimensional linear subspace onto which the projection of the centered random vector retains as much of the original information as possible. 
Throughout this paper, we assume that $1 \leq r < d$. 
Under an arbitrary distribution $\dP$, PCA can be cast into the following optimization model,
\begin{equation}
	\label{opt:pca}
	\min_{X \in \Odr} 
	\hkh{
		\sfR_{\dP} (X)
		:= 
		\dE_{\dP}
		\fkh{ \norm{ \dkh{ I_d - X X\zz } \dkh{ \xi - \dE_{\dP} \fkh{\xi} } }_2^2 }
	},
\end{equation}
where $\Odr = \hkh{X \in \Rdr \mid X\zz X = I_r}$ denotes the Stiefel manifold \cite{Absil2008Optimization,Boumal2023Introduction,Wang2022Decentralized} in $\Rdr$. 
For a candidate subspace spanned by the columns of $X \in \Odr$, $\sfR_{\dP} (X)$ quantifies the reconstruction error incurred by dimension reduction under the distribution $\dP$.

Replacing $\dP$ in problem~\eqref{opt:pca} with $\dP\uast$ yields the true PCA model that we seek to solve. 
In most practical scenarios, the true distribution $\dP\uast$ is unknown and can not be captured precisely, leaving us without the essential information required to solve it exactly. 
A standard remedy is to substitute $\dP\uast$ with the empirical distribution $\dP_n\ucirc = n^{- 1} \sum_{i = 1}^n \delta_{\xi_i\ucirc}$ constructed from the observed samples $\{\xi_i\ucirc\}_{i \in [n]}$ with $[n] = \hkh{1, \ldots, n}$, where $\delta_{\xi_i\ucirc}$ denotes the Dirac distribution concentrated at $\xi_i\ucirc$. 
To account explicitly for the distributional uncertainty arising from this approximation, we consider an infinity-Wasserstein distributionally robust formulation of PCA. 
Below is the definition of the infinity-Wasserstein distance.

\begin{definition}[\cite{Bertsimas2023Data,Givens1984Class}]
	Let $\sH (\dP_1, \dP_2)$ be the set of all joint distributions on $\Rd \times \Rd$ with first and second marginals $\dP_1$ and $\dP_2$ on $\Rd$. 
	The infinity-Wasserstein distance between $\dP_1$ and $\dP_2$ is defined as
	\begin{equation*}
		\dW_\infty (\dP_1, \dP_2)
		:= \inf_{\dJ \in \sH (\dP_1, \dP_2)}
		\esssup_{(x, y) \sim \dJ} \norm{x - y}_2,
	\end{equation*}
	where $\esssup_{(x, y) \sim \dJ}$ denotes the essential supremum with respect to $\dJ$, namely the smallest extended-real number $t$ for which $\norm{x - y}_2 \leq t$ holds $\dJ$-almost surely. 
\end{definition}

Given a radius $\rho \geq 0$, the corresponding infinity-Wasserstein ambiguity set centered at the empirical distribution is
\begin{equation*}
	\sB_\infty (\dP_n\ucirc, \rho)
	:= \hkh{\dP \mid \dW_\infty (\dP_n\ucirc, \dP) \leq \rho}.
\end{equation*}
Then we consider the worst-case reconstruction error over the ambiguity set and obtain the following DRO model,
\begin{equation}
	\label{opt:infwpca}
	\min_{X \in \Odr}
	\hkh{
		\sfV_{n, \rho} (X)
		:= \sup_{\dP \in \sB_\infty (\dP_n\ucirc, \rho)} 
		\sfR_{\dP} (X)
	},
\end{equation}
which is called {\it infinity-Wasserstein PCA} and hereafter referred to as InfWPCA. 
A notable feature of problem~\eqref{opt:infwpca} is that its dependence on the underlying probability distribution is fundamentally different from that in conventional Wasserstein DRO models. 
A standard formulation typically involves an objective of the form $\dE_{\dP} [\sfU (\xi)]$, where the loss function $\sfU$ is fixed independently of $\dP$. 
The expectation is therefore a linear functional of the probability distribution. 
In the PCA setting considered here, however, the distribution enters the objective twice. 
In addition to determining the outer expectation, it also determines the mean $\dE_{\dP} [\xi]$ appearing inside $\sfR_{\dP}$ defined by \eqref{opt:pca}. 
The reconstruction error exhibits a quadratic dependence on the distribution through the interaction between these two expectations. 
Consequently, the inner worst-case problem in \eqref{opt:infwpca} is different from the standard Wasserstein DRO setting, which calls for a more careful reformulation.

\subsection{Exact Finite-Dimensional Reformulation}

\label{subsec:projected_cen_repre}

The nonlinear dependence on $\dP$ identified in the preceding subsection results from the distribution-dependent mean $\mathbb{E}_{\dP}[\xi]$.
We observe that $\mathbb{E}_{\dP}[\xi]$ can be separated from the reconstruction error through a variational representation. 
This structure allows us to evaluate the Wasserstein ambiguity and derive an exact finite-dimensional reformulation of problem~\eqref{opt:infwpca}. 

We introduce a center variable $z \in \Rd$ in place of the distribution-dependent mean $\mathbb{E}_{\dP}[\xi]$ in the reconstruction error. 
The fixed-center loss is written as
\begin{equation*}
	\sfL_X (\xi, z) = \norm{\dkh{I_d - X X\zz} \dkh{\xi - z}}_2^2.
\end{equation*}
Then it follows from the standard variance decomposition that
\begin{equation}
	\label{eq:variance_identity}
    \dE_{\dP} \fkh{\sfL_X (\xi, z)}
    = \sfR_{\dP} (X)
    + \norm{
    \dkh{I_d - X X\zz}
    \dkh{z - \dE_{\dP} \fkh{\xi}}
    }_2^2,
\end{equation}
which further yields the variational representation
\begin{equation}
	\label{eq:variational}
    \sfR_{\dP} (X)
    = \inf_{z \in \Rd} \dE_{\dP} \fkh{\sfL_X (\xi, z)}.
\end{equation}
As a result, the reconstruction error can be expressed as the minimum of an expected loss whose integrand no longer depends on $\dP$. 
Based on this representation, we now establish the desired exact reformulation.

\begin{theorem}
	\label{thm:exact_reformulation}
	For every $X \in \Odr$ and $\rho \geq 0$, it holds that
	\begin{equation}
		\label{eq:exact_reformulation}
		\sfV_{n, \rho} (X)
		= \min_{z \in \Rd} \hkh{
			\sfF_{n, \rho} (X, z)
			:= \dfrac{1}{n} \sum_{i = 1}^n
			\dkh{\norm{\sfr_i (X, z)}_2 + \rho}^2
		},
	\end{equation}
	where $\sfr_i (X, z) = \dkh{I_d - X X\zz} \dkh{\xi_i\ucirc - z}$. 
\end{theorem}

\begin{proof}
	If $\rho = 0$, the ambiguity set $\sB_\infty (\dP_n\ucirc, 0)$ collapses to the singleton $\{\dP_n\ucirc\}$. 
	Then the relationship \eqref{eq:exact_reformulation} holds trivially. 
	Next, our attention is restricted to the case $\rho > 0$. 
	It is clear that $\hat{z} = \dE_{\dP} [\xi]$ is a global minimizer of the optimization problem in \eqref{eq:variational}. 
	Let $\mu_n\ucirc = n^{-1} \sum_{i = 1}^n \xi_i\ucirc$ be the empirical mean. 
	For any $\dP \in \sB_\infty (\dP_n\ucirc, \rho)$, we have
	\begin{equation*}
		\norm{ \hat{z} - \mu_n\ucirc }_2
		= \norm{ \dE_{\dP} \fkh{\xi} - \dE_{\dP_n\ucirc} \fkh{\xi} }_2
		\leq \inf_{\dJ \in \sH (\dP_n\ucirc, \dP)} \dE_{(\xi\ucirc, \xi) \sim \dJ} \fkh{\norm{\xi - \xi\ucirc}_2}
		\leq \dW_\infty (\dP_n\ucirc, \dP)
		\leq \rho.
	\end{equation*}
	Thus, the center variable $z$ can be restricted to a compact set $\cB (\mu_n\ucirc, \rho) = \hkh{z \in \Rd \mid \norm{z - \mu_n\ucirc}_2 \leq \rho}$. 
	Moreover, it holds that
	\begin{equation*}
		\sfR_{\dP} (X)
		= \min_{z \in \Rd} \dE_{\dP} \fkh{\sfL_X (\xi, z)}
		= \min_{z \in \cB (\mu_n\ucirc, \rho)} \dE_{\dP} \fkh{\sfL_X (\xi, z)}. 
	\end{equation*}
	Since $z \mapsto \sfL (\xi, z)$ is continuous and convex for each $\xi$, the compactness and regularity conditions required by \cite[Theorem~3(II)]{Zhang2025Short} are satisfied. 
	Based on the definition of $\dP_n\ucirc$, we then arrive at
	\begin{equation}
		\label{eq:dual}
		\sfV_{n, \rho} (X)
		= \min_{z \in \Rd} \frac{1}{n} \sum_{i = 1}^n
		\sup_{\norm{\delta_i}_2 \leq \rho}
		\sfL_X (\xi_i\ucirc + \delta_i, z).
	\end{equation}
	For a fixed pair $(X, z) \in \Odr \times \Rd$, the triangle inequality yields that
	\begin{equation*}
		\begin{aligned}
			\sfL_X (\xi_i\ucirc + \delta_i, z)
			= {} & \norm{ \dkh{I_d - X X\zz} \dkh{\xi_i\ucirc + \delta_i - z} }_2^2 \\
			\leq {} & \dkh{ \norm{\sfr_i (X, z)}_2 + \norm{\dkh{I_d - X X\zz} \delta_i}_2 }^2
			\leq \dkh{ \norm{\sfr_i (X, z)}_2 + \rho }^2,
		\end{aligned}
	\end{equation*}
	where the last inequality follows from the orthogonality of $X$ and $\norm{\delta_i}_2 \leq \rho$. 
	Now we choose
	\begin{equation*}
		\delta_i\ustar =
		\left\{
		\begin{aligned}
			& \rho \sfr_i (X, z) / \norm{\sfr_i (X, z)}_2,
			&& \mbox{if~} \sfr_i (X, z) \neq 0, \\
			& \rho (I_d - X X\zz) v / \norm{(I_d - X X\zz) v}_2,
			&& \mbox{if~} \sfr_i (X, z) \neq 0, \\
		\end{aligned}
		\right.
	\end{equation*}
	where $v \in \Rd$. 
	Then it can be readily verified that $\norm{\delta_i\ustar}_2 \leq \rho$ and
	\begin{equation*}
		\begin{aligned}
			\sfL_X (\xi_i\ucirc + \delta_i\ustar, z)
			= \norm{ \dkh{I_d - X X\zz} \dkh{\xi_i\ucirc + \delta_i\ustar - z} }_2^2
			= \dkh{ \norm{\sfr_i (X, z)}_2 + \rho }^2,
		\end{aligned}
	\end{equation*}
	which indicates that $\sup_{\norm{\delta_i}_2 \leq \rho} \sfL_X (\xi_i\ucirc + \delta_i, z) = \dkh{ \norm{\sfr_i (X, z)}_2 + \rho }^2$ for $i \in [n]$. 
	Combining this equality with the relationship \eqref{eq:dual}, we finally obtain the desired result \eqref{eq:exact_reformulation}. 
	The proof is completed. 
\end{proof}

A key implication of Theorem~\ref{thm:exact_reformulation} is that the original infinite-dimensional DRO problem admits an exact finite-dimensional reformulation. 
The equivalent optimization model takes the form
\begin{equation}
    \label{opt:infwpca-eq}
    \min_{X \in \Odr, z \in \Rd} \sfF_{n, \rho} (X, z),
\end{equation}
where $\sfF_{n, \rho}$ is defined in \eqref{eq:exact_reformulation}. 
This representation combines the empirical reconstruction error with an average residual-norm term arising from the infinity-Wasserstein geometry. 

Alongside the subspace variable $X \in \Odr$, the exact reformulation in \eqref{opt:infwpca-eq} introduces an auxiliary center variable $z \in \Rd$. 
The next result characterizes the behavior of its optimal solutions. 

\begin{theorem}
	\label{thm:center_attainment}
	Let $\cC_X := \conv \hkh{\dkh{I_d - X X\zz} \xi_i\ucirc \mid i \in [n]}$, where $\conv\, \cC$ denotes the convex hull of a set $\cC$. 
	For any fixed $X \in \Odr$, the minimization of $\sfF_{n, \rho} (X, z)$ with respect to $z \in \cS_X^\perp := \hkh{u \in \Rd \mid X\zz u = 0}$ has a unique solution $z\ustar_X \in \cC_X$. 
	Moreover, there exists a uniform constant $\nu := \max_{i \in [n]} \norm{\xi_i\ucirc}_2$ such that $\norm{z\ustar_X}_2 \leq \nu$ for all $X \in \Odr$. 
	The complete set of global solutions to the minimization of $\sfF_{n, \rho} (X, z)$ over $z \in \Rd$ is $z\ustar_X + \cS_X$ with $\cS_X := \hkh{X v \mid v \in \Rr}$. 
\end{theorem}

\begin{proof}
	For any $z \in \cS_X^\perp$, we have $(I_d - X X\zz) z = z$. 
	Expanding the restriction of $\sfF_{n, \rho} (X, z)$ to $z \in \cS_X^\perp$ gives
	\begin{equation*}
		\frac{1}{n} \sum_{i = 1}^n
		\norm{\dkh{I_d - X X\zz}\xi_i\ucirc - z}_2^2
		+ \frac{2\rho}{n} \sum_{i = 1}^n
		\norm{\dkh{I_d - X X\zz}\xi_i\ucirc - z}_2
		+ \rho^2.
	\end{equation*}
	The first term is strongly convex and coercive on $\cS_X^\perp$, while the second term is convex on $\cS_X^\perp$. 
	Therefore, the minimization of $\sfF_{n, \rho} (X, z)$ with respect to $z \in \cS_X^\perp := \hkh{u \in \Rd \mid X\zz u = 0}$ has a unique solution, denoted by $z\ustar_X$.

	Let $\proj_{\cC_X} (z)$ be the Euclidean projection of $z \in \cS_X^\perp$ onto $\cC_X$. 
	It is clear that $\proj_{\cC_X} (z) \in \cS_X^\perp$. 
	For each $y \in \cC_X$, the projection inequality yields
	\begin{equation}
		\label{eq:projection}
		\norm{y - z}_2^2
		\geq \norm{y - \proj_{\cC_X} (z)}_2^2 + \norm{z - \proj_{\cC_X} (z)}_2^2.
	\end{equation}
	We assume, on the contrary, that $z\ustar_X \notin \cC_X$. 
	Then $z\ustar_X \neq \proj_{\cC_X} (z\ustar_X)$. 
	It follows from \eqref{eq:projection} that $\norm{\dkh{I_d - X X\zz}\xi_i\ucirc - z\ustar_X}_2 > \norm{\dkh{I_d - X X\zz}\xi_i\ucirc - \proj_{\cC_X} (z\ustar_X)}_2$ for each $i \in [n]$, which further implies that $\sfF_{n, \rho} (X, z\ustar_X) > \sfF_{n, \rho} (X, \proj_{\cC_X} (z\ustar_X))$. 
	This leads to a contradiction to the global optimality of $z\ustar_X$. 
	As a result, we have $z\ustar_X \in \cC_X$. 
	The definitions of $\nu$ and $\cC_X$ together guarantee that $\norm{z\ustar_X}_2 \leq \nu$.

	Finally, we choose an arbitrary global solution $\bar{z}$ to the minimization of $\sfF_{n, \rho} (X, z)$ over $z \in \Rd$. 
	Obviously, it holds that $(I_d - X X\zz) \bar{z} \in \cS_X^\perp$
	\begin{equation*}
		\sfF_{n, \rho} (X, (I_d - X X\zz) \bar{z}) 
		= \sfF_{n, \rho} (X, \bar{z})
		\leq \sfF_{n, \rho} (X, z),
	\end{equation*}
	for all $z \in \cS_X^\perp$. 
	Hence, $(I_d - X X\zz) \bar{z}$ is a global solution to the minimization of $\sfF_{n, \rho} (X, z)$ with respect to $z \in \cS_X^\perp$. 
	The definition and the uniqueness of $z\ustar_X$ implies that $(I_d - X X\zz) \bar{z} = z\ustar_X$. 
	It follows that $\bar{z} = (I_d - X X\zz) \bar{z} + X X\zz \bar{z} \in z\ustar_X + \cS_X$. 
	Conversely, it can be readily verified that any point in $z\ustar_X + \cS_X$ is a global solution $\bar{z}$ to the minimization of $\sfF_{n, \rho} (X, z)$ over $z \in \Rd$. 
	We complete the proof. 
\end{proof}

Theorem~\ref{thm:center_attainment} demonstrates that, in problem~\eqref{opt:infwpca-eq}, the minimization with respect to $z$ can be restricted to the compact set $\hkh{z \in \Rd \mid \norm{z}_2 \leq \nu}$ without altering the optimal value. 
Together with the compactness of $\Odr$ and the continuity of $\sfF_{n, \rho}$, there always exists a global minimizer for problem~\eqref{opt:infwpca-eq}. 

We henceforth refer to $z\ustar_X$ as the optimal center, which can be interpreted as a generalized sample mean induced by the infinity-Wasserstein geometry. 
This center also plays a pivotal role in the explicit construction of a worst-case distribution and in the subsequent convergence analysis of the proposed algorithm.

\subsection{Worst-Case Distribution }

\label{subsec:worst_distribution}

After establishing the exact reformulation, we now turn to the distribution that attains the worst-case reconstruction error in the DRO model \eqref{opt:infwpca}. 
The main difficulty lies in reconciling the extremal perturbations with the requirement that their resulting mean be consistent with the optimal center characterized by Theorem~\ref{thm:center_attainment}. 

Given the optimal center $z\ustar_X$, we define $\cI_0 = \hkh{i \in [n] \mid \sfr_i (X, z\ustar_X) = 0}$. 
According to the optimality condition for the minimization of $\sfF_{n, \rho} (X, z)$ with respect to $z$, there exist vectors $u_i \in \cS_X^\perp$ for $i \in \cI_0$ satisfying $\norm{u_i}_2 \leq 1$ and
\begin{equation}
	\label{eq:subgradient}
	\sum_{i \notin \cI_0}
	\dkh{
		\sfr_i (X, z\ustar_X)
		+ \frac{\rho \sfr_i (X, z\ustar_X)}{\norm{\sfr_i (X, z\ustar_X)}_2}
	}
	+ \rho \sum_{i \in \cI_0} u_i
	= 0.
\end{equation}
The following theorem gives a constructive worst-case distribution.

\begin{theorem}
	\label{thm:inner_worst_case_attainment}
	For each $i \in \cI_0$, we set $q_i = u_i / \norm{u_i}_2$ when $u_i \neq 0$ and choose any unit vector $q_i \in \cS_X^\perp$ when $u_i = 0$, where $u_i \in \cS_X^\perp$ satisfying $\norm{u_i}_2 \leq 1$ and \eqref{eq:subgradient}. 
	Let $\alpha_i = \dkh{1 + \norm{u_i}_2} / 2$ and
	\begin{equation*}
		\dQ_i\udiamond = 
		\left\{
		\begin{aligned}
			& \delta_{\xi_i\ucirc + \rho \sfr_i (X, z\ustar_X) / \|\sfr_i (X, z\ustar_X)\|_2},
			&& \mbox{if~} i \notin \cI_0, \\
			& \alpha_i \delta_{\xi_i\ucirc + \rho q_i}
			+ \dkh{1 - \alpha_i} \delta_{\xi_i\ucirc - \rho q_i},
			&& \mbox{if~} i \in \cI_0.
		\end{aligned}
		\right.
	\end{equation*}
	Then the distribution $\dP_X\udiamond = n^{-1} \sum_{i = 1}^n \dQ_i\udiamond$ belongs to $\sB_\infty (\dP_n\ucirc, \rho)$ and satisfies $\sfV_{n, \rho} (X) = \sfR_{\dP_X\udiamond} (X)$ for any $X \in \Odr$ and $\rho \geq 0$. 
\end{theorem}

\begin{proof}
	If $\rho = 0$, the constructed distribution $\dP_X\udiamond$ reduces to $\dP_n\ucirc$. 
	The assertions are immediate in this case. 
	Suppose now that $\rho > 0$. 
	According to the definitions of $\alpha_i$ and $q_i$, we have
	\begin{equation*}
		\dE_{\dQ_i\udiamond} \fkh{\xi} 
		= \xi_i\ucirc + (2 \alpha_i - 1) \rho q_i
		= \xi_i\ucirc + \rho \norm{u_i}_2 q_i
		= \xi_i\ucirc + \rho u_i,
	\end{equation*}
	for each $i \in \cI_0$. 
	Recall that $\mu_n\ucirc = n^{-1} \sum_{i = 1}^n \xi_i\ucirc$ is the empirical mean. 
	Then it can be readily verified that
	\begin{equation*}
		\begin{aligned}
			\dE_{\dP_X\udiamond} \fkh{\xi}
			= {} & \dfrac{1}{n} \sum_{i \notin \cI_0} \dkh{ \xi_i\ucirc + \dfrac{\rho \sfr_i (X, z\ustar_X)}{\norm{\sfr_i (X, z\ustar_X)}_2} }
			+ \dfrac{1}{n} \sum_{i \in \cI_0} \dkh{ \xi_i\ucirc + \rho u_i } \\
			= {} & \mu_n\ucirc + \dfrac{1}{n} \sum_{i \notin \cI_0} \dfrac{\rho \sfr_i (X, z\ustar_X)}{\norm{\sfr_i (X, z\ustar_X)}_2}
			+ \dfrac{\rho}{n} \sum_{i \in \cI_0} u_i.
		\end{aligned}
	\end{equation*}
	As a direct consequence of \eqref{eq:subgradient}, we proceed to show that
	\begin{equation*}
		\dE_{\dP_X\udiamond} \fkh{\xi}
		= \mu_n\ucirc
		- \dfrac{1}{n} \sum_{i \notin \cI_0} \sfr_i (X, z\ustar_X)
		= \mu_n\ucirc
		- \dfrac{1}{n} \sum_{i = 1}^{n} \sfr_i (X, z\ustar_X),
	\end{equation*}
	where the last equality follows from the definition of $\cI_0$. 
	Substituting the expression $\sfr_i (X, z\ustar_X) = \dkh{I_d - X X\zz} \dkh{\xi_i\ucirc - z\ustar_X}$ into the above relationship yields that
	\begin{equation}
		\label{eq:center}
		\dE_{\dP_X\udiamond} \fkh{\xi}
		= \dkh{I_d - X X\zz} z\ustar_X + X X\zz \mu_n\ucirc
		= z\ustar_X + X X\zz \mu_n\ucirc.
	\end{equation}
	For each $i \in \cI_0$, since $q_i \in \cS_X^\perp$ satisfies $\norm{q_i}_2 = 1$ and $\sfr_i (X, z\ustar_X) = 0$, we have
	\begin{equation*}
		\begin{aligned}
			\dE_{\dQ_i\udiamond} \fkh{ \norm{ \dkh{ I_d - X X\zz } \dkh{ \xi - z\ustar_X } }_2^2 }
			= {} & \alpha_i \norm{ \sfr_i (X, z\ustar_X) + \rho q_i }_2^2 \\
			& + \dkh{1 - \alpha_i} \norm{ \sfr_i (X, z\ustar_X) - \rho q_i }_2^2 \\
			= {} & \alpha_i \norm{ \rho q_i}_2^2
			+ \dkh{1 - \alpha_i} \norm{ \rho q_i }_2^2 \\
			= {} & \rho^2
			= \dkh{\sfr_i (X, z\ustar_X) + \rho}^2.
		\end{aligned}
	\end{equation*}
	In addition, it holds for all $i \notin \cI_0$ that
	\begin{equation*}
		\begin{aligned}
			\dE_{\dQ_i\udiamond} \fkh{ \norm{ \dkh{ I_d - X X\zz } \dkh{ \xi - z\ustar_X } }_2^2 }
			= {} & \norm{ \sfr_i (X, z\ustar_X) + \dfrac{\rho \sfr_i (X, z\ustar_X)}{\norm{\sfr_i (X, z\ustar_X)}_2} }_2^2 \\
			= {} & \dkh{\norm{\sfr_i (X, z\ustar_X)}_2 + \rho}^2.
		\end{aligned}
	\end{equation*}
	Then it follows from the relationship \eqref{eq:center} and the definition of $\dP_X\udiamond$ that
	\begin{equation*}
		\begin{aligned}
			\sfR_{\dP_X\udiamond} (X)
			= {} & \dE_{\dP_X\udiamond} \fkh{ \norm{ \dkh{ I_d - X X\zz } \dkh{ \xi - \dE_{\dP_X\udiamond} \fkh{\xi} } }_2^2 } \\
			= {} & \dE_{\dP_X\udiamond} \fkh{ \norm{ \dkh{ I_d - X X\zz } \dkh{ \xi - z\ustar_X } }_2^2 } \\
			= {} & \dfrac{1}{n} \sum_{i = 1}^{n} \dkh{\norm{\sfr_i (X, z\ustar_X)}_2 + \rho}^2
			= \sfF_{n, \rho} (X, z\ustar_X).
		\end{aligned}
	\end{equation*}
	According to Theorem~\ref{thm:center_attainment}, we know that $z\ustar_X$ is a global minimizer of $\sfF_{n, \rho} (X, z)$ over $z \in \Rd$. 
	Hence, Theorem~\ref{thm:exact_reformulation} guarantees that $\sfR_{\dP_X\udiamond} (X) = \sfF_{n, \rho} (X, z\ustar_X) = \min_{z \in \Rd} \sfF_{n, \rho} (X, z) = \sfV_{n, \rho} (X)$. 
	
	It remains to verify that $\dP_X\udiamond \in \sB_\infty (\dP_n\ucirc, \rho)$. 
	We consider the joint distribution $\dJ\udiamond = n^{-1} \sum_{i = 1}^n \delta_{\xi_i\ucirc} \otimes \dQ_i\udiamond$, which has marginal distributions $\dP_n\ucirc$ and $\dP_X\udiamond$. 
	Moreover, a straightforward verification reveals that
	\begin{equation*}
		\dJ\udiamond \dkh{ \hkh{ (x,y) \mid \norm{x - y}_2 > \rho } }
		= \dfrac{1}{n} \sum_{i = 1}^n
		\dQ_i\udiamond
		\dkh{ \hkh{ y \mid \norm{\xi_i\ucirc - y}_2 > \rho } }
		= 0,
	\end{equation*}
	which is followed by
	\begin{equation*}
		\dW_\infty (\dP_n\ucirc, \dP_X\udiamond)
		\leq \esssup_{(x, y) \sim \dJ\udiamond} \norm{x - y}_2
		\leq \rho.
	\end{equation*}
	Therefore, we have $\dP_X\udiamond \in \sB_\infty (\dP_n\ucirc, \rho)$ and complete the proof. 
\end{proof}

Theorem~\ref{thm:inner_worst_case_attainment} furnishes an explicit construction of a distribution within the ambiguity set that attains the worst-case reconstruction error. 
To the best of our knowledge, no comparable constructive result has been reported in the existing literature on DRO formulations of PCA. 
This observation sets the stage for a detailed comparison in the next subsection.

\subsection{Comparison with Existing Models}

\label{subsec:comparison_type2}

The purpose of this subsection is to discuss the distinctions between the proposed infinity-Wasserstein model and existing order-$2$ Wasserstein models of PCA. 
The central question is whether the choice of Wasserstein geometry can generate a nontrivial robustification effect without sacrificing an exact finite-dimensional reformulation. 

We begin by revisiting the order-$2$ Wasserstein DRO model. 
Let $G \in \mathbb{R}^{d \times d}$ be a positive definite matrix. 
The weighted order-$2$ Wasserstein distance \cite{Givens1984Class} is given by
\begin{equation*}
    \dW_2^G (\dP_1, \dP_2)
    = \dkh{
    \inf_{\dJ \in \sH (\dP_1, \dP_2)}
    \dE_{(x, y) \sim \dJ}
    \fkh{
    \dkh{x - y}\zz
    G
    \dkh{x - y}
    }
    }^{1 / 2}.
\end{equation*}
The associated worst-case reconstruction error takes the form
\begin{equation}
	\label{eq:order2wdro}
	\sfT_{n, \rho}^G (X)
	= \sup_{\dP \in \sB_2^G (\dP_n\ucirc, \rho)} \sfR_{\dP} (X),
\end{equation}
with the ambiguity set $\sB_2^G (\dP_n\ucirc, \rho) = \hkh{\dP \mid \dW_2^G (\dP_n\ucirc, \dP) \leq \rho}$ and the radius $\rho \geq 0$.

When $G = I_d$, the superscript $G$ is omitted from $\dW_2^G$, $\sB_2^G$, and $\sfT_{n, \rho}^G$. 
For this case, Wang et al.~\cite{Wang2025Enhancing} establish the exact identity
\begin{equation*}
    \sfT_{n, \rho} (X)
    = \dkh{
    \sqrt{\tr \dkh{\dkh{I_d - X X\zz}\Sigma_n\ucirc}}
    + \rho
    }^2,
\end{equation*}
where $\Sigma_n\ucirc = n^{-1} \sum_{i = 1}^n \dkh{\xi_i\ucirc - \mu_n\ucirc} \dkh{\xi_i\ucirc - \mu_n\ucirc}\zz$ is the empirical covariance matrix. 
Since $t \mapsto \dkh{\sqrt{t} + \rho}^2$ is strictly increasing on $[0, \infty)$, we have
\begin{equation*}
    \argmin_{X \in \Odr} \; \sfT_{n, \rho} (X)
    = \argmin_{X \in \Odr} \; \tr \dkh{\dkh{I_d - X X\zz} \Sigma_n\ucirc}.
\end{equation*}
Thus, the isotropic order-$2$ robustification leaves the principal subspaces unchanged. 

For a general positive definite matrix $G$, an explicit expression for the worst-case reconstruction error defined in \eqref{eq:order2wdro} is typically unavailable. 
To circumvent this difficulty, Xu et al.~\cite{Xu2026Distributionally} propose the following upper-bound surrogate,
\begin{equation}
	\label{eq:w2_upper_bound}
    \sfJ_{n, \rho}^G (X)
    = \sqrt{\tr \dkh{\dkh{I_d - X X\zz} \Sigma_n\ucirc}}
    + \rho
    \sqrt{
    \lambda_{\max}
    \dkh{
    G^{-1 / 2}
    \dkh{I_d - X X\zz}
    G^{-1 / 2}
    }
    },
\end{equation}
where $\lambda_{\max} (\cdot)$ represents the largest eigenvalue of a given matrix. 
It has been proved that $\sfT_{n, \rho}^G (X) \leq (\sfJ_{n, \rho}^G (X))^2$ for all $X \in \Odr$. 
Although this weighted order-$2$ construction can select a principal subspace different from that of nominal PCA, the minimization of its surrogate remains challenging due to the nonlinear coupling between the empirical reconstruction error and the largest-eigenvalue penalty.

Unlike existing order-$2$ Wasserstein models, the InfWPCA model \eqref{opt:infwpca} reconciles the exact finite-dimensional reformulation \eqref{opt:infwpca-eq} with the capacity to yield a robust principal subspace. 

\begin{example}
	\label{ex:winfty_differs_nominal}
	We consider $d = 2$, $r = 1$, $n = 4$, and $\rho = 1$.
	The samples are given by $\xi_1\ucirc = (1, 3 / 2)$, $\xi_2\ucirc = (1, - 3 / 2)$, and $\xi_3\ucirc = \xi_4\ucirc = (- 1, 0)$. 
	A direct calculation gives $\mu_n\ucirc = 0$ and $\Sigma_n\ucirc = \Diag \dkh{1, 9 / 8}$. 
	Let $e_i$ be the $i$-th standard basis vector in $\Rd$. 
	It is clear that $e_2$ is the eigenvector corresponding to the largest eigenvalue of $\Sigma_n\ucirc$, which is a global minimizer of nominal PCA. 
	Based on the relationship \eqref{eq:exact_reformulation}, we have
	\begin{equation*}
		\sfV_{n, \rho} (e_2)
		= \frac{1}{4} \sum_{i = 1}^4 \dkh{1 + \rho}^2
		= 4.
	\end{equation*}
	and
	\begin{equation*}
		\sfV_{n, \rho} (e_1)
		= \frac{1}{4}
		\fkh{
			2 \dkh{\frac{3}{2} + \rho}^2
			+ 2 \rho^2
		}
		= \frac{29}{8}
		< 4
		= \sfV_{n, \rho} (e_2).
	\end{equation*}
	Consequently, the nominal principal subspace is not optimal for the InfWPCA model. 
\end{example}

Both nominal PCA and the order-$2$ Wasserstein surrogate in \eqref{eq:w2_upper_bound} for a fixed matrix $G$ depend on the empirical distribution solely through its covariance matrix. 
By contrast, the InfWPCA model retains information beyond the first two empirical moments, as illustrated by the following example.

\begin{example}
	\label{ex:equal_moments}
	Let $d = 2$, $r = 1$, and $n = 4$. 
	We consider two sets of samples, including $\hkh{e_2, e_2, - e_2, - e_2}$ and $\hkh{\sqrt{2} e_2, - \sqrt{2} e_2, 0, 0}$. 
	The samples in both sets have zero mean and covariance matrix $\Diag (0, 1)$. 
	By symmetry and Theorem~\ref{thm:center_attainment}, the optimal center at $X = e_1$ is zero in both cases. 
	Their worst-case reconstruction errors are $1 + 2 \rho + \rho^2$ and $1 + \sqrt{2} \rho + \rho^2$, respectively. 
	Hence, the InfWPCA model is not determined by the first two empirical moments whenever $\rho > 0$. 
\end{example}

\section{Finite-Sample Guarantee}

\label{sec:statistical_properties}

When the true distribution has unbounded support, it necessarily lies outside every finite-radius infinity-Wasserstein ambiguity set centered at the empirical distribution. 
Nevertheless, this lack of distributional containment does not preclude meaningful finite-sample guarantees of the DRO model~\eqref{opt:infwpca}. 
By exploiting the exact reformulation derived in Section~\ref{sec:model_analysis}, we establish a high-probability and uniform upper bound on the reconstruction error under the true distribution, directly in terms of the distributionally robust objective function.



Throughout this section, the observed sample $\xi_1\ucirc, \dotsc, \xi_n\ucirc$ are independent and identically distributed according to the true distribution $\dP\uast$ with the following assumption.

\begin{assumption}
	\label{ass:stat_bounded_kurtosis}
	There exists a finite constant $\kappa > 0$ such that $\dE_{\dP\uast} [ \|\xi - \dE_{\dP\uast} [\xi]\|_2^4 ] \leq \kappa$. 
\end{assumption}

Assumption~\ref{ass:stat_bounded_kurtosis} imposes a finite fourth-moment condition and accommodates a broad class of probability distributions, including sub-Gaussian distributions and certain heavy-tailed distributions. 
Notably, it permits polynomially decaying tails, as exemplified by multivariate Student distributions with more than four degrees of freedom. 
This requirement is considerably weaker than the distributional assumptions employed in existing finite-sample results. 
Nguyen et al.~\cite{Nguyen2020Distributionally} demonstrate that a true distribution supported on a bounded domain lies within an infinity-Wasserstein ambiguity set with high probability. 
When the true distribution is sub-Gaussian, Yue et al.~\cite{Yue2026Geometric} establish analogous guarantees for distributionally robust covariance estimators with divergence-based ambiguity sets. 
Rather than relying on distributional containment, the finite-sample analysis developed here is closer in spirit to that of Gao~\cite{Gao2023Finite}, which derives nonasymptotic performance guarantees for order-$1$ and order-$2$ Wasserstein DRO under a transportation-information inequality. 
Our results extend this performance-oriented perspective to the infinity-Wasserstein model of PCA in substantially less restrictive distributional settings.

Unless otherwise specified, all expectations and probabilities are taken with respect to the true distribution $\dP\uast$ in this section, together with its induced distributions on transformed random variables and independent samples. 
Expectations involving auxiliary random variables are understood to be taken jointly over these variables and the observed samples.

\subsection{Whitening and Augmented Residual Representation}

\label{subsec:whitening}

To facilitate the finite-sample analysis, we construct a unified representation of reconstruction residuals through covariance-support whitening and deterministic-coordinate augmentation. 
This construction yields the moment properties required for the uniform probabilistic bounds.

Let $\mu\uast$ and $\Sigma\uast$ be the mean vector and the covariance matrix of the true distribution $\dP\uast$, respectively. 
If $\Sigma\uast = 0$, all the observed samples coincide with $\mu\uast$ almost surely and $\sfR_{\dP\uast} (X) = 0$ for every $X \in \Odr$. 
The desired guarantee is therefore immediate in this case. 
We henceforth restrict attention to $\Sigma\uast \neq 0$. 
Let $s \geq 1$ be the rank of $\Sigma\uast$. 
The singular value decomposition of $\Sigma\uast$ can be expressed as $\Sigma\uast = U\uast \Lambda\uast (U\uast)\zz$, where $U\uast \in \Ods$ is an orthogonal matrix and $\Lambda\uast \in \Rss$ is a diagonal matrix whose diagonal entries are arranged in descending order as $\lambda\uast_1 \geq \dotsb \geq \lambda\uast_s > 0$. 
In addition, the random vector $\xi - \mu\uast$ belongs to the range space of $U\uast$ almost surely, since its projection onto the orthogonal complement has zero second moment. 
Accordingly, there exists a random vector $\varsigma \in \Rs$ such that $\xi - \mu\uast = U\uast (\Lambda\uast)^{1/2} \varsigma$. 
It can be readily verified that
\begin{equation*}
	\varsigma = (\Lambda\uast)^{-1/2} (U\uast)\zz \dkh{\xi - \mu\uast},
\end{equation*}
which satisfies $\dE [\varsigma] = 0$ and $\dE [\varsigma \varsigma\zz] = I_s$. 
By Assumption~\ref{ass:stat_bounded_kurtosis}, we have
\begin{equation}
	\label{eq:stat_whitening}
	\dE \fkh{\abs{u\zz \varsigma}^4} 
	\leq \dE \fkh{ \norm{(\Lambda\uast)^{-1/2} u}_2^4   \norm{(U\uast)\zz \dkh{\xi - \mu\uast}}_2^4}
	\leq \dfrac{\kappa}{(\lambda\uast_s)^2} \norm{u}_2^4
	= \bar{\kappa} \norm{u}_2^4,
\end{equation}
for all $u \in \Rs$. 
Here, $\bar{\kappa} := \kappa / (\lambda\uast_s)^2$ is a finite constant. 
A direct application of Jensen's inequality gives $\kappa \geq (\dE_{\dP\uast} [ \|\xi - \mu\uast\|_2^2 ])^2 = (\tr (\Sigma\uast))^2 \geq (\lambda\uast_s)^2$, from which $\bar{\kappa} \geq 1$ immediately follows.

In the next step, we adjoin a deterministic coordinate to $\varsigma$ by setting $\zeta = (\varsigma\zz, 1)\zz \in \bR^{s+1}$. 
A straightforward verification reveals that
\begin{equation}
	\label{eq:stat_augmented}
	\dE\fkh{\zeta \zeta\zz} = I_{s + 1}. 
\end{equation}
For $v = (u\zz, a)\zz \in \bR^{s + 1}$, the elementary inequality $|x + y|^4 \leq 8 |x|^4 + 8 |y|^4$ and the relationship \eqref{eq:stat_whitening} yield
\begin{equation}
	\label{eq:stat_augmented_fourth}
	\dE \fkh{\abs{v\zz \zeta}^4}
	\leq 8 \bar{\kappa} \norm{u}_2^4 + 8 a^4
	\leq 8 \bar{\kappa} \dkh{\norm{u}_2^2 + a^2}^2
	= 8 \bar{\kappa} \norm{v}_2^4.
\end{equation}
Thus, the augmented vector is isotropic in the second-moment sense and retains uniform fourth-moment control.

Let $\varsigma_i$ and $\zeta_i$ be the samples of $\varsigma$ and $\zeta$ corresponding to $\xi_i\ucirc$, respectively. 
Given $X\in\Odr$ and $z\in\Rd$, we define
\begin{equation*}
	B_{X, z} := \fkh{ \dkh{I_d - X X\zz} U\uast (\Lambda\uast)^{1/2} \quad - \dkh{I_d - X X\zz} \dkh{z - \mu\uast} } \in \bR^{d \times (s + 1)}.
\end{equation*}
Then it holds that $B_{X, z} \zeta_i = \sfr_i (X, z)$. 
Moreover, the augmentation procedure produces the exact identity
\begin{equation}
	\label{eq:stat_augmented_identities}
	\norm{B_{X, z}}\fs = \norm{ \dkh{I_d - X X\zz} U\uast (\Lambda\uast)^{1/2} }\fs + \norm{ \dkh{I_d - X X\zz} \dkh{z - \mu\uast} }_2^2.
\end{equation}
These relationship establish an essential link between empirical residual moments and reconstruction error under the true distribution.

\subsection{Lower Bounds for Empirical Residual Moments}

Building on the augmented moment representation, we derive two complementary high-probability lower bounds for empirical reconstruction residuals. 
The first uses a one-sided quadratic-form inequality to control the squared residuals, while the second employs a truncation-based argument to establish a uniform lower bound on their average norms. 
Together, these results provide the probabilistic foundation for the finite-sample guarantee.


We first establish a high-probability lower bound for the empirical second-moment matrix of the augmented samples in the Loewner order. 
Drawing on the one-sided quadratic-form inequality of Oliveira~\cite{Oliveira2016Lower}, this result requires only uniform fourth-moment control.

\begin{lemma}
	\label{lem:stat_moment_event}
	Let $\zeta_1, \dotsc, \zeta_n$ be independent samples of a random vector $\zeta \in \bR^{s + 1}$ satisfying \eqref{eq:stat_augmented} and \eqref{eq:stat_augmented_fourth}. 
	For any $t > 0$, it holds with probability at least $1 - e^{- t}$ that
	\begin{equation}
		\label{eq:stat_matrix_lower_event}
		M_n := \dfrac{1}{n} \sum_{i = 1}^n \zeta_i \zeta_i\zz \succeq \max \hkh{1 - \beta_{n} (t), 0} I_{s + 1},
	\end{equation}
	where
	\begin{equation*}
		\beta_{n} (t) := 18 \sqrt{\frac{2 \bar{\kappa} (s + 1 + 2 t + 2 \log 2)}{n}},
	\end{equation*}
	and the notation $A \succeq B$ represents that $A - B$ is positive semidefinite. 
\end{lemma}

\begin{proof}
	For convenience, we define $H_i := \zeta_i \zeta_i\zz \succeq 0$. 
	Then it follows from \eqref{eq:stat_augmented} that $\dE [H_i] = I_{s + 1}$. 
	Let $\zeta_{i, j}$ be the $j$-th entry of $\zeta_i$. 
	The entries of $H_i$ have finite second moments, since the Cauchy-Schwarz inequality and \eqref{eq:stat_augmented_fourth} imply $\dE [\zeta_{i, j}^2 \zeta_{i, k}^2] \leq 8 \bar{\kappa}$ for every $j$ and $k$. 
	Moreover, for all $v \in \bR^{s + 1}$, we have
	\begin{equation*}
		\sqrt{ \dE \fkh{\dkh{v\zz H_i v}^2} }
		= \sqrt{ \dE\fkh{\abs{v\zz \zeta_i}^4} }
		\leq \sqrt{8 \bar{\kappa}} \norm{v}_2^2
		= \sqrt{8 \bar{\kappa}} v\zz \dE \fkh{H_i} v.
	\end{equation*}
	As a result, Oliveira's one-sided quadratic-form inequality~\cite[Theorem~4.1]{Oliveira2016Lower} applies with moment parameter $h = \sqrt{8 \bar{\kappa}} > 1$ and confidence level $\delta = e^{-t}$. 
	The desired result \eqref{eq:stat_matrix_lower_event} follows immediately, and we complete the proof. 
\end{proof}



Lemma~\ref{lem:stat_moment_event} provides a lower bound on the empirical squared residuals, but the resulting estimate may fall short of the reconstruction error under the true distribution. 
To bridge this gap, we proceed to establish a lower bound on the empirical residual norms. 

\begin{lemma}
	\label{lem:stat_pair_variation}
	Suppose that $\zeta_1, \dotsc, \zeta_n$ satisfy the assumptions of Lemma~\ref{lem:stat_moment_event}. 
	Let
	\begin{equation}
		\label{eq:stat_pair_sample_size}
		n \geq \max \hkh{ 8192 \bar{\kappa} \dkh{s + 1}, 512 \bar{\kappa}^2 t },
	\end{equation}
	for all $t \geq 1$. 
	Then, with probability at least $1 - e^{- t}$, we have
	\begin{equation}
		\label{eq:stat_matrix_first_moment}
		\dfrac{1}{n} \sum_{i = 1}^n \norm{B \zeta_i}_2
		\geq \dfrac{1}{8 \sqrt{2 \bar{\kappa}}} \norm{B}\ff,
	\end{equation}
	for every $B \in \bR^{d \times (s + 1)}$. 
\end{lemma}

\begin{proof}
	Let $v \in \bR^{s + 1}$ be a unit vector and $T = v\zz \zeta$. 
	The moment assumptions give $\dE [T^2] = 1$ and $\dE [T^4] \leq 8 \bar{\kappa}$. 
	Interpolation between the first, second, and fourth moments yields
	\begin{equation*}
		\dE \fkh{\abs{T}}
		\geq \dfrac{ \dkh{ \dE \fkh{T^2} }^{3 / 2} }{ \dkh{\dE \fkh{T^4}}^{1/2} }
		\geq \dfrac{1}{2 \sqrt{2 \bar{\kappa}}}.
	\end{equation*}
	We define a bounded and Lipschitz continuous function $\varphi (x) = \min \{\abs{x}, \sqrt{2 \bar{\kappa}}\}$. 
	The elementary inequality $\max \{\abs{x} - \sqrt{2 \bar{\kappa}}, 0\} \leq x^2 / (4 \sqrt{2 \bar{\kappa}})$ for $x \in \bR$ shows that
	\begin{equation}
		\label{eq:stat_clipped_population}
		\dE \fkh{\varphi (T)}
		= \dE \fkh{\abs{T}} - \dE \fkh{\max\{\abs{T} - \sqrt{2 \bar{\kappa}}, 0\}}
		\geq \frac{1}{2 \sqrt{2 \bar{\kappa}}} - \frac{1}{4 \sqrt{2 \bar{\kappa}}}
		= \dfrac{1}{4 \sqrt{2 \bar{\kappa}}}.
	\end{equation}
	To control all directions at once, we set
	\begin{equation*}
		D_n = \max \hkh{ 0, \sup_{\norm{v}_2 = 1} \dkh{ \dE \fkh{\varphi \dkh{v\zz \zeta}} - \dfrac{1}{n} \sum_{i = 1}^n \varphi \dkh{v\zz \zeta_i} } }.
	\end{equation*}
	The supremum is measurable because the empirical paths and expectations are continuous in $v$, allowing restriction to a countable dense subset of the unit sphere. 
	The standard symmetrization inequality and Rademacher contraction principle~\cite[Theorem~A.3]{Dalalyan2018Optimal}, together with $\varphi (0) = 0$ and the unit Lipschitz constant of $\varphi$, give
	\begin{equation*}
		\begin{aligned}
			\dE \fkh{D_n}
			\leq {} & \dfrac{2}{n} \dE \fkh{ \sup_{\norm{v}_2 = 1} \abs{\sum_{i = 1}^n \varepsilon_i \varphi \dkh{v\zz \zeta_i}} }
			\leq \dfrac{4}{n} \dE \fkh{ \sup_{\norm{v}_2 = 1} \abs{v\zz \sum_{i = 1}^n \varepsilon_i \zeta_i} } \\
			= {} & \dfrac{4}{n} \dE \fkh{ \norm{\sum_{i = 1}^n \varepsilon_i \zeta_i}_2 }
			\leq 4 \sqrt{\dfrac{s + 1}{n}},
		\end{aligned}
	\end{equation*}
	where $\varepsilon_1, \dotsc, \varepsilon_n$ are independent Rademacher signs, independent of the samples. 
	The last step follows from Jensen's inequality, sign orthogonality, and $\dE [\norm{\zeta}_2^2] = s + 1$.

	Since $0 \leq \varphi (x) \leq \sqrt{2 \bar{\kappa}}$ for all $x \in \bR$, replacing any sample changes $D_n$ by at most $\sqrt{2 \bar{\kappa}} / n$. 
	According to McDiarmid's bounded-differences inequality~\cite[Theorem~A.2]{Dalalyan2018Optimal}, we have
	\begin{equation*}
		D_n \leq \dE \fkh{D_n} + \sqrt{\dfrac{\bar{\kappa} t}{n}},
	\end{equation*}
	with probability at least $1 - e^{- t}$.
	The two requirements in \eqref{eq:stat_pair_sample_size} guarantee
	\begin{equation*}
		4 \sqrt{\dfrac{s + 1}{n}} \leq \dfrac{1}{16 \sqrt{2 \bar{\kappa}}},
		\quad\mbox{and}\quad
		\sqrt{\dfrac{\bar{\kappa} t}{n}} \leq \frac{1}{16 \sqrt{2 \bar{\kappa}}},
	\end{equation*}
	respectively.
	Thus, on the same event, one has $D_n \leq 1 / (8 \sqrt{2 \bar{\kappa}})$. Combining this bound with \eqref{eq:stat_clipped_population} and $\abs{x} \geq \varphi (x)$ establishes
	\begin{equation}
		\label{eq:stat_pair_variation}
		\dfrac{1}{n} \sum_{i = 1}^n \abs{v\zz \zeta_i}
		\geq \dfrac{1}{8 \sqrt{2 \bar{\kappa}}},
	\end{equation}
	whenever $\norm{v}_2 = 1$.

	For the desired result \eqref{eq:stat_matrix_first_moment}, the case $B = 0$ is immediate. 
	Otherwise, we consider the eigenvalue decomposition of $B\zz B$ as
	\begin{equation*}
		B\zz B = \sum_{j = 1}^{s + 1} \lambda_j u_j u_j\zz,
	\end{equation*}
	where $\lambda_1, \dotsc, \lambda_{s + 1}$ are eigenvalues and $u_1, \dotsc, u_{s + 1}$ are orthonormal eigenvectors. 
	It follows from the Cauchy-Schwarz inequality that
	\begin{equation*}
		\sum_{j = 1}^{s + 1} \lambda_j \abs{u_j \zz w}
		\leq \dkh{\sum_{j = 1}^{s + 1} \lambda_j}^{1/2} \dkh{\sum_{j = 1}^{s + 1} \lambda_j \dkh{u_j\zz w}^2}^{1 / 2}
		= \norm{B}\ff \norm{B w}_2,
	\end{equation*}
	for all $w \in \bR^{s + 1}$. 
	Applying \eqref{eq:stat_pair_variation} to each $u_j$ on the same probability event yields
	\begin{equation*}
		\dfrac{1}{n} \sum_{i = 1}^{n} \norm{B \zeta_i}_2
		\geq \frac{1}{\norm{B}\ff} \sum_{j = 1}^{s + 1} \lambda_j \dkh{ \dfrac{1}{n} \sum_{i = 1}^{n} \abs{u_j\zz \zeta_i} }
		\geq \dfrac{1}{8 \sqrt{2 \bar{\kappa}} \norm{B}\ff}
		\sum_{j = 1}^{s + 1} \lambda_j
		= \dfrac{\norm{B}\ff}{8 \sqrt{2 \bar{\kappa}}}.
	\end{equation*}
	The proof is completed.
\end{proof}

In the subsequent analysis, Lemma~\ref{lem:stat_pair_variation} will be instrumental in showing that the residual-norm term induced by the infinity-Wasserstein formulation compensates for potential shortfalls in the empirical squared residuals. 
This result plays a crucial role in establishing the desired finite-sample performance guarantee.

\subsection{Uniform Out-of-Sample Certificate}
\label{subsec:stat_rootn}

The preceding results reveal how the residual structure of the InfWPCA model can be exploited to control the reconstruction error under the true distribution. 
We now translate these insights into a uniform finite-sample performance guarantee.

\begin{theorem}
	\label{thm:stat_performance}
	Suppose that Assumption~\ref{ass:stat_bounded_kurtosis} holds. 
	For every integer $n \geq 1$ and every $t \geq 1$, we define
	\begin{equation}
		\label{eq:stat_radius_remainder}
		\rho_{n} (t) := \dfrac{144 \kappa}{(\lambda\uast_s)^2} \sqrt{\dfrac{ \lambda\uast_1 \bar{s} \dkh{s + 1 + 2 t + 2 \log 2} }{n}},
	\end{equation}
	where $\bar{s} := \min \{s, d - r\}$. 
	Then, with probability at least $1 - 2 e^{- t}$, it holds that
	\begin{equation}
		\label{eq:stat_main_certificate}
		\sfR_{\dP\uast} (X) \leq \sfV_{n, \rho} (X),
	\end{equation}
	for all $X \in \Odr$ and $\rho \geq \rho_{n} (t)$. 
\end{theorem}

\begin{proof}
	For every $X \in \Odr$, the matrix $(I_d - X X\zz) U\uast (\Lambda\uast)^{1/2}$ has rank at most $\bar{s}$, and its largest singular value does not exceed $(\lambda\uast_1)^{1/2}$. 
	Then it can be readily verified that
	\begin{equation}
		\label{eq:stat_rank_bound}
		\sfR_{\dP\uast} (X)
		= \tr \dkh{(I_d - X X\zz) \Sigma\uast}
		= \norm{(I_d - X X\zz) U\uast (\Lambda\uast)^{1/2}}\fs 
		\leq \lambda\uast_1 \bar{s}.
	\end{equation}
	We first assume that $144 \bar{\kappa} \sqrt{(s + 1 + 2 t + 2 \log 2) / n} \geq 1$. 
	Then it holds that $\rho_{n} (t) \geq \sqrt{\lambda\uast_1 \bar{s}}$.
	According to Theorem~\ref{thm:exact_reformulation} and the estimate \eqref{eq:stat_rank_bound}, we arrive at
	\begin{equation*}
		\sfV_{n, \rho} (X)
		\geq \rho^2 
		\geq \lambda\uast_1 \bar{s}
		\geq \sfR_{\dP\uast} (X),
	\end{equation*}
	for $X \in \Odr$ and $\rho \geq \rho_{n} (t)$. 
	
	It remains to consider the case $144 \bar{\kappa} \sqrt{(s + 1 + 2 t + 2 \log 2) / n} < 1$. 
	Then we have
	\begin{equation*}
		n > 20736 \bar{\kappa}^2 \dkh{s + 1 + 2 t + 2 \log 2} \geq \max \hkh{8192 \bar{\kappa} \dkh{s + 1}, 512 \bar{\kappa}^2 t},
	\end{equation*}
	where the latter inequality uses $\bar{\kappa} \geq 1$. 
	Hence, the sample-size condition \eqref{eq:stat_pair_sample_size} of Lemma~\ref{lem:stat_pair_variation} is automatically satisfied. 
	Under Assumption~\ref{ass:stat_bounded_kurtosis}, the random vector $\zeta$ constructed in Section~\ref{subsec:whitening} adheres to the moment conditions \eqref{eq:stat_augmented} and \eqref{eq:stat_augmented_fourth}. 
	All the requirements of Lemma~\ref{lem:stat_moment_event} and Lemma~\ref{lem:stat_pair_variation} are fulfilled. 
	Let $\cE_1$ and $\cE_2$ denote the events in \eqref{eq:stat_matrix_lower_event} and \eqref{eq:stat_matrix_first_moment}, respectively. 
	Then the event $\cE_1 \cap \cE_2$ occurs with probability at least $1 - 2 e^{- t}$. 
	We fix a realization in this intersection. 
	Recall from Section~\ref{subsec:whitening} that $\sfr_i (X, z) = B_{X, z} \zeta_i$. 
	It follows from Lemma~\ref{lem:stat_moment_event} that
	\begin{equation*}
		\dfrac{1}{n} \sum_{i = 1}^{n} \norm{\sfr_i (X, z)}_2^2
		= \dfrac{1}{n} \sum_{i = 1}^{n} \norm{B_{X, z} \zeta_i}_2^2
		= \tr \dkh{ B_{X, z} M_n B_{X, z}\zz }
		\geq \max \hkh{1 - \beta_{n} (t), 0} \norm{B_{X, z}}\fs.
	\end{equation*}
	In light of Lemma~\ref{lem:stat_pair_variation}, we further obtain that
	\begin{equation*}
		\dfrac{1}{n} \sum_{i = 1}^{n} \norm{\sfr_i (X, z)}_2
		= \dfrac{1}{n} \sum_{i = 1}^{n} \norm{B_{X, z} \zeta_i}_2
		\geq \dfrac{1}{8 \sqrt{2 \bar{\kappa}}} \norm{B_{X, z}}\ff.
	\end{equation*}
	The definitions of $\beta_{n} (t)$ and $\rho_{n} (t)$ imply that $\rho_{n} (t) = 4 \sqrt{2 \bar{\kappa} \bar{s} \lambda\uast_1} \beta_{n} (t)$. 
	For every $\rho \geq \rho_{n} (t)$, the estimate \eqref{eq:stat_rank_bound} guarantees that
	\begin{equation}
		\label{eq:linear_residual}
		\dfrac{\rho}{4 \sqrt{2 \bar{\kappa}}} \sqrt{\sfR_{\dP\uast} (X)}
		\geq \beta_{n} (t) \sqrt{\lambda\uast_1 \bar{s} \sfR_{\dP\uast} (X)}
		\geq \beta_{n} (t) \sfR_{\dP\uast} (X).
	\end{equation}
	Expanding the function $\sfF_{n, \rho}$ in \eqref{eq:exact_reformulation} gives
	\begin{equation*}
		\begin{aligned}
			\sfF_{n, \rho} (X, z)
			={} & \dfrac{1}{n} \sum_{i = 1}^{n} \norm{\sfr_i (X, z)}_2^2 + \dfrac{2 \rho}{n} \sum_{i = 1}^{n} \norm{\sfr_i (X, z)}_2 + \rho^2 \\
			\geq {} & \max \hkh{1 - \beta_{n} (t), 0} \norm{B_{X, z}}\fs + \dfrac{\rho}{4 \sqrt{2 \bar{\kappa}}} \norm{B_{X,z}}\ff + \rho^2.
		\end{aligned}
	\end{equation*}
	For arbitrary $X\in\Odr$ and $z\in\Rd$, the identity \eqref{eq:stat_augmented_identities} and the relationship \eqref{eq:stat_rank_bound} yield $\norm{B_{X, z}} \fs \geq \sfR_{\dP\uast} (X)$. 
	Consequently, we have
	\begin{equation*}
		\begin{aligned}
			\sfF_{n, \rho} (X, z)
			\geq {} & \max \hkh{1 - \beta_{n} (t), 0} \sfR_{\dP\uast} (X) + \dfrac{\rho}{4 \sqrt{2 \bar{\kappa}}} \sqrt{\sfR_{\dP\uast} (X)} + \rho^2 \\
			\geq {} & \dkh{ \max \hkh{1 - \beta_{n} (t), 0} + \beta_{n} (t) } \sfR_{\dP\uast} (X) + \rho^2
			\geq \sfR_{\dP\uast} (X),
		\end{aligned}
	\end{equation*}
	where the second inequality follows from \eqref{eq:linear_residual}. 
	Taking the infimum over $z$ and invoking Theorem~\ref{thm:exact_reformulation} establishes \eqref{eq:stat_main_certificate}, which completes the proof. 
\end{proof}

Theorem~\ref{thm:stat_performance} demonstrates that the objective function of the InfWPCA model \eqref{opt:infwpca} provides a uniform upper bound on the reconstruction error under the true distribution with high probability, thereby delivering a rigorous finite-sample guarantee for the generalization performance. 
In its proof, both probabilistic bounds hold simultaneously over all candidate subspaces and centers on a common event. 
Therefore, the out-of-sample certificate \eqref{eq:stat_main_certificate} remains valid for data-dependent choices of $X$ and $\rho$. 

\section{Algorithm Design}

\label{sec:algorithm}

The exact finite-dimensional reformulation established in Section~\ref{sec:model_analysis} provides a tractable representation of the proposed InfWPCA model, although the resulting optimization problem~\eqref{opt:infwpca-eq} remains nonsmooth and nonconvex. 
Such problems can be addressed by general-purpose Riemannian optimization algorithms, including subgradient methods \cite{Li2021Weakly} and smoothing approaches \cite{Beck2023Dynamic,Wang2026Adaptive,Zhang2024Riemannian}. 
The particular structure of samplewise residuals, however, naturally lends itself to an iteratively reweighted majorization-minimization scheme, enabling a highly efficient algorithm tailored to problem~\eqref{opt:infwpca-eq}. 
At each iteration, we devise a safeguarded quadratic majorant, whose global minimization reduces to computing a weighted sample mean and performing a weighted PCA update. 
Beyond its computational simplicity, this construction also lays the foundation for the convergence analysis. 
Throughout this section, the sample size $n \geq 1$ and the radius $\rho > 0$ are fixed.


\subsection{A Safeguarded Quadratic Majorant}

\label{subsec:mm_majorization}

We begin with a smooth quadratic majorant for the residual-norm term in problem~\eqref{opt:infwpca-eq}. 
The safeguard introduced in this procedure prevents singular reweighting when a residual vanishes.

For every $t \geq 0$ and every $\gamma > 0$, we have $2 t \leq t^2 / \gamma + \gamma$. 
Accordingly, each residual-norm term admits the following upper bound,
\begin{equation*}
	\dkh{\norm{\sfr_i (X, z)}_2 + \rho}^2
	\leq \dkh{1 + \dfrac{\rho}{\gamma_i}} \norm{\sfr_i (X, z)}_2^2
    + \rho \gamma_i
    + \rho^2,
\end{equation*}
for all $\gamma_i > 0$. 
Let $(X^{(k)}, z^{(k)}) \in \Odr \times \Rd$ be the current iterate at iteration $k$. 
For each $i \in [n]$, we introduce the reweighting coefficient $\gamma_i^{(k)}$ by
\begin{equation}
	\label{eq:gamma}
	\gamma_i^{(k)}
	:= \max \hkh{
		\norm{\sfr_i (X^{(k)}, z^{(k)})}_2,
		\tau^{(k)}
	},
\end{equation}
where $\tau^{(k)} > 0$ is the safeguard parameter to be specified later. 
This mechanism guarantees that every weight is finite and positive. 
Then the quadratic majorant is defined by
\begin{equation}
	\label{eq:iteration_majorant_model}
	\sfQ^{(k)} (X, z)
	:= \dfrac{1}{n} \sum_{i = 1}^n
	\dkh{
		\omega_i^{(k)} \norm{\sfr_i (X, z)}_2^2
		+ \rho \gamma_i^{(k)}
		+ \rho^2
	},
\end{equation}
with
\begin{equation}
	\label{eq:weight}
	\omega_i^{(k)}
	:= 1 + \dfrac{\rho}{\gamma_i^{(k)}},
\end{equation}
for each $i \in [n]$. 
The resulting reweighting scheme places greater emphasis on smaller residuals, as dictated by the structure of the infinity-Wasserstein formulation. 
A straightforward verification reveals that $\sfF_{n, \rho} (X, z) \leq \sfQ^{(k)} (X, z)$ on $\Odr \times \Rd$. 

For the convergence analysis, it is convenient to identify the merit function to which the quadratic majorant is tangent. 
For $\tau > 0$ and $t \geq 0$, we define the smoothing envelope by
\begin{equation*}
	\sfN_\tau (t)
	:= \min_{\gamma \geq \tau}
	\hkh{
		\dfrac{t^2}{2\gamma}
		+ \dfrac{\gamma}{2}
	}
	= 
	\begin{cases}
		\dfrac{t^2}{2\tau} + \dfrac{\tau}{2},
		& \mbox{if~} 0 \leq t < \tau, \\
		t,
		& \mbox{if~} t \geq \tau,
	\end{cases}
\end{equation*}
where the minimum is uniquely attained at $\gamma = \max \{t, \tau\}$. 
It can be readily verified that $0 \leq \sfN_\tau (t) - t \leq \tau / 2$ for all $t \geq 0$. 
This envelope gives rise to the following merit function,
\begin{equation}
	\label{eq:merit}
	\sfE_\tau (X, z)
	:= \dfrac{1}{n} \sum_{i = 1}^n
	\dkh{
		\norm{\sfr_i (X, z)}_2^2
		+ 2 \rho \sfN_\tau \dkh{\norm{\sfr_i (X, z)}_2}
		+ \rho^2
	},
\end{equation}
which satisfies $\sfF_{n, \rho} (X, z) \leq \sfE_\tau (X, z) \leq \sfF_{n, \rho} (X, z) + \rho\tau$ for all $(X, z) \in \Odr \times \Rd$. 
Moreover, the choice of $\gamma_i^{(k)}$ in \eqref{eq:gamma} attains the minimum in the smoothing-envelope representation at the current residual norm. 
Combining this fact with \eqref{eq:iteration_majorant_model} yields
\begin{equation}
	\label{eq:merit_majorization_by_Q}
	\sfE_{\tau^{(k)}} (X, z)
	\leq \sfQ^{(k)} (X, z),
	\quad\mbox{and}\quad
	\sfE_{\tau^{(k)}} (X^{(k)}, z^{(k)})
	= \sfQ^{(k)} (X^{(k)}, z^{(k)}).
\end{equation}
Thus, the majorant $\sfQ^{(k)}$ serves as a global upper bound of the merit function and is tangent to it at the current iterate.

\subsection{Closed-Form Reweighted Updates}

\label{subsec:closed_form_mm_updates}

With the quadratic majorant in hand, we now turn to its minimization over $\Odr \times \Rd$. 
Despite the joint dependence on the subspace and center variables, a global minimizer admits a particularly simple construction.

We construct a global minimizer of $\sfQ^{(k)}$ in two successive steps. 
The weighted sample mean is given by
\begin{equation}
	\label{eq:z_update}
	z^{(k + 1)}
	:=
	\dfrac{
		\sum_{i = 1}^n \omega_i^{(k)} \xi_i\ucirc
	}{
		\sum_{i = 1}^n \omega_i^{(k)}
	}.
\end{equation}
Since every weight is positive, $z^{(k + 1)}$ belongs to the convex hull of the observed samples. 
Then we introduce the weighted covariance matrix as
\begin{equation}
	\label{eq:weighted_covariance}
	\Sigma_{\omega}^{(k)}
	:= \dfrac{1}{n} \sum_{i = 1}^n
	\omega_i^{(k)}
	\dkh{\xi_i\ucirc - z^{(k + 1)}}
	\dkh{\xi_i\ucirc - z^{(k + 1)}}\zz.
\end{equation}
With the center variable $z$ fixed at $z^{(k + 1)}$, the minimization of $\sfQ^{(k)}$ with respect to $X$ can be recast as the following weighted PCA problem,
\begin{equation}
	\label{eq:X_update_eig}
	X^{(k + 1)}
	\in \argmax_{X \in \Odr}
	\; \tr \dkh{X\zz \Sigma_{\omega}^{(k)} X}.
\end{equation}
Hence, $X^{(k + 1)}$ can be chosen to form an orthonormal basis of an $r$-dimensional invariant subspace associated with the $r$ largest eigenvalues of $\Sigma_{\omega}^{(k)}$. 
The proposition below corroborates that $(X^{(k + 1)}, z^{(k + 1)})$ is a global minimizer of $\sfQ^{(k)}$.

\begin{proposition}
	\label{prop:joint_majorant_minimization}
	For the constructions in \eqref{eq:z_update} and \eqref{eq:X_update_eig}, we have
	\begin{equation*}
		\dkh{X^{(k + 1)}, z^{(k + 1)}}
		\in
		\argmin_{X \in \Odr, z \in \Rd}
		\sfQ^{(k)} (X, z).
	\end{equation*}
\end{proposition}

\begin{proof}
	We first fix an arbitrary $X \in \Odr$. 
	It follows from the update scheme \eqref{eq:z_update} of $z^{(k + 1)}$ that
	\begin{equation*}
		\sum_{i = 1}^n
		\omega_i^{(k)}
		\dkh{\xi_i\ucirc - z^{(k + 1)}}
		= 0.
	\end{equation*}
	Applying the weighted variance identity then yields
	\begin{equation*}
		\sum_{i = 1}^n
		\omega_i^{(k)}
		\norm{\sfr_i (X, z)}_2^2
		= \sum_{i = 1}^n
		\omega_i^{(k)}
		\norm{\sfr_i (X, z^{(k + 1)})}_2^2
		+ \sum_{i = 1}^n \omega_i^{(k)}
		\norm{
			\dkh{I_d - X X\zz}
			\dkh{z - z^{(k + 1)}}
		}_2^2,
	\end{equation*}
	for all $z \in \Rd$. 
	Since the second term is always nonnegative, $z^{(k + 1)}$ is a global minimizer of $\sfQ^{(k)} (X, z)$ with respect to $z$ for every fixed $X$. 
	
	Now we set $z = z^{(k + 1)}$. 
	The terms $\rho\gamma_i^{(k)} + \rho^2$ in \eqref{eq:iteration_majorant_model} are independent of $X$, and the remaining weighted residual sum can be written as
	\begin{equation*}
		\sum_{i = 1}^n \omega_i^{(k)} \norm{\sfr_i (X, z^{(k + 1)})}_2^2
		= \sum_{i = 1}^n \omega_i^{(k)} \norm{\xi_i\ucirc - z^{(k + 1)}}_2^2
		- n \tr \dkh{X\zz \Sigma_{\omega}^{(k)} X}.
	\end{equation*}
	Therefore, the minimization of $\sfQ^{(k)} (X, z^{(k + 1)})$ with respect to $X$ is equivalent to the weighted PCA problem in \eqref{eq:X_update_eig}. 
	Combining the optimality of \eqref{eq:z_update} for every fixed $X$ with the optimality of \eqref{eq:X_update_eig} gives $\sfQ^{(k)} (X^{(k + 1)}, z^{(k + 1)}) \leq \sfQ^{(k)} (X, z)$ for all $(X, z) \in \Odr \times \Rd$. 
	The proof is completed.
\end{proof}

The closed-form reweighted updates in \eqref{eq:z_update} and \eqref{eq:X_update_eig} form the computational core of an efficient algorithm for solving problem~\eqref{opt:infwpca-eq}. 
The remaining task is to specify a safeguard schedule and design a certificate-based output selection rule.

\subsection{Complete Framework}

\label{subsec:wellposed_mm}

We conclude this section by presenting the complete framework of the proposed algorithm, which is termed {\it Iteratively Reweighted Wasserstein PCA} and abbreviated as IRWPCA.

At iteration $k$, the safeguard parameter follows the polynomial decay schedule $\tau^{(k)} = \tau^{(0)} (k + 1)^{- \sigma}$ with $\tau^{(0)} > 0$ and $\sigma \in (0, 1)$. 
This choice prevents singular reweighting at vanishing residuals while progressively reducing the effect of safeguard. 
As will be shown later, the decay exponent $\sigma$ determines the order of iteration complexity. 
The next iterates $z^{(k + 1)}$ and $X^{(k + 1)}$ are then updated by the reweighted steps \eqref{eq:z_update} and \eqref{eq:X_update_eig}, respectively. 
Moreover, IRWPCA maintains a certificate-based output sequence $\{\bar{X}^{(k)}\}$, whose selection is guided by the following majorization gap,
\begin{equation}
	\label{eq:majorization_gap}
	\Delta^{(k)}
	:= \sfQ^{(k)} (X^{(k + 1)}, z^{(k + 1)})
	- \sfE_{\tau^{(k)}} (X^{(k + 1)}, z^{(k + 1)}).
\end{equation}
It follows from \eqref{eq:merit_majorization_by_Q} that $\Delta^{(k)} \geq 0$. 
After each iteration, the output is chosen from the latter half of the iteration history by minimizing the normalized gap $\Delta^{(j)} / \tau^{(j)}$ over $j \in \{\lfloor (k + 1) / 2 \rfloor, \dotsc, k\}$. 
This quantity controls the squared stationarity residual through an explicit upper bound, which will be elaborated upon in the next section. 
Consequently, the resulting selection criterion favors an iterate with a stronger stationarity certificate. 
Algorithm~\ref{alg:irwpca} summarizes the whole procedure of IRWPCA.

\begin{algorithm2e}[ht]
	\caption{Iteratively Reweighted Wasserstein PCA (IRWPCA).}
	\label{alg:irwpca}
	
	\KwIn{$X^{(0)} \in \Odr$, $z^{(0)} \in \Rd$, $\tau^{(0)} > 0$, $\sigma \in (0, 1)$, and $K \geq 1$.}
	
	\For{$k = 0, 1, 2, \dotsc, K - 1$}{
		
		Choose $\tau^{(k)} = \tau^{(0)} (k + 1)^{- \sigma}$.
		
		Set $\gamma_i^{(k)}$ and $\omega_i^{(k)}$ by \eqref{eq:gamma} and \eqref{eq:weight} for each $i \in [n]$, respectively. 
		
		Compute $z^{(k + 1)}$ according to \eqref{eq:z_update}.
		
		Form $\Sigma_{\omega}^{(k)}$ based on \eqref{eq:weighted_covariance} and update $X^{(k + 1)}$ from \eqref{eq:X_update_eig}. 
		
		Calculate $\Delta^{(k)}$ using \eqref{eq:majorization_gap}. 
		
		Select $\sfm (k) \in \argmin_{j \in \{\lfloor (k + 1) / 2 \rfloor, \dotsc, k\}} \hkh{ \Delta^{(j)} / \tau^{(j)} }$. 
		
		Assign $\bar{X}^{(k + 1)} = X^{(\sfm (k) + 1)}$, $\bar{z}^{(k + 1)} = z^{(\sfm (k) + 1)}$, and $\bar{\tau}^{(k)} = \tau^{(\sfm (k))}$. 
		
	}
	
	\KwOut{$(\bar{X}^{(K)}, \bar{z}^{(K)})$.}
\end{algorithm2e}

\section{Convergence Analysis}

\label{sec:convergence}

This section delves into a rigorous convergence analysis of Algorithm~\ref{alg:irwpca}, with the sample size $n \geq 1$ and the radius $\rho > 0$ held fixed. 
Although iteratively reweighted PCA methods have been extensively studied for other types of problems \cite{Lerman2025Global,Lerman2018Fast,Neumayer2020Robust,Park2016Iteratively}, existing convergence results often rely on specific structural assumptions on the data or concern only smoothed surrogate objectives. 
Within the framework of Riemannian nonsmooth analysis, we establish the Clarke stationarity of accumulation points and derive the corresponding iteration complexity, without imposing additional nondegeneracy or distributional assumptions on the data. 
Notably, while IRWPCA operates on the joint variables $(X, z) \in \Odr \times \Rd$ in problem~\eqref{opt:infwpca-eq}, our theoretical guarantees pertain to the original objective function $\sfV_{n, \rho}$ of the InfWPCA model~\eqref{opt:infwpca}.

\subsection{Riemannian Clarke Subdifferential}

\label{subsec:value_subdifferential}

The first step is to characterize the Riemannian Clarke subdifferential \cite{Hosseini2018Line,Li2021Weakly,Yang2014Optimality} of $\sfV_{n,\rho}$ based on the exact representation in \eqref{eq:exact_reformulation}, which is a natural extension of the Clarke subdifferential \cite{Clarke1990Optimization,Rockafellar1998Variational} in the Euclidean space. 
For $X \in \Odr$, the tangent space \cite{Absil2008Optimization,Boumal2023Introduction} of the Stiefel manifold is $\cT_X \Odr := \{ H \in \Rdr \mid X\zz H + H\zz X = 0 \}$. 
The Riemannian metric on $\cT_X \Odr$ is induced from the Frobenius inner product. 
The orthogonal projection operator onto $\cT_X \Odr$ is $\proj_{\cT_X \Odr} (H) := H - X (X\zz H + H\zz X) / 2$. 
We adopt the following definition of the Riemannian Clarke subdifferential.

\begin{definition}
	Suppose that $\sff: \Rdr \to \bR$ is a locally Lipschitz continuous function. 
	Let $\Omega \lrmR (\sff) := \{ X \in \Odr \mid \mbox{$\sff$ is differentiable at $X$} \}$. 
	Then the Riemannian Clarke subdifferential of $\sff$ at $X \in \Odr$ is defined as
	\begin{equation*}
		\partial\lrmR \sff (X) = \conv \hkh{H \in \cT_{X} \Odr \mid \grad\, \sff (X_t) \to H, \Omega\lrmR (\sff) \ni X_t \to X},
	\end{equation*}
	where $\grad\, \sff (X) := \proj_{\cT_X \Odr} (\nabla \sff (X))$ represents the Riemannian gradient of $\sff$ at a differentiable point $X \in \Odr$. 
\end{definition}

Let $\partial\, \sff$ be the Clarke subdifferential of $\sff$ in the Euclidean space. 
When $\sff$ is Clarke regular \cite{Yang2014Optimality} at $X \in \Odr$, its Riemannian Clarke subdifferential admits the representation $\partial\lrmR \sff (X) = \proj_{\cT_X \Odr} (\partial\, \sff (X))$. 
Recall from Theorem~\ref{thm:center_attainment} that $\cC_X = \conv \{\dkh{I_d - X X\zz} \xi_i\ucirc \mid i \in [n]\}$ and $\cS_X^\perp = \{u \in \Rd \mid X\zz u = 0\}$. 
For $X \in \Odr$, $z \in \cC_X$, and $\tau \geq 0$, we introduce the following admissible direction set,
\begin{equation}
	\label{eq:direction}
	\cU_\tau (X, z)
	:= \hkh{ u = (u_1, \ldots, u_n) \;\middle|\;
		\begin{aligned}
			& u_i = \dfrac{\sfr_i (X, z)}{\norm{\sfr_i (X, z)}_2}, \mbox{ if } \norm{\sfr_i (X, z)}_2 > \tau \\
			& u_i \in \cS_X^\perp, \norm{u_i}_2 \leq 1, \mbox{ if } \norm{\sfr_i (X, z)}_2 \leq \tau
		\end{aligned}
	}.
\end{equation}
For $u \in \cU_\tau (X, z)$, the corresponding first-order components are given by
\begin{equation}
	\label{eq:value_first_order_vectors}
	\left\{
	\begin{aligned}
		\sfG_X (X, z, u)
		:= {} & - \dfrac{2}{n} \sum_{i = 1}^n \dkh{\sfr_i (X, z) + \rho u_i} \dkh{\xi_i\ucirc}\zz X, \\
		\sfG_z (X, z, u)
		:= {} & - \dfrac{2}{n} \sum_{i = 1}^n \dkh{\sfr_i (X, z) + \rho u_i}.
	\end{aligned}
	\right.
\end{equation}
Since $z \in \cC_X \subseteq \cS_X^\perp$ and $u_i \in \cS_X^\perp$, the matrix $\sfG_X (X, z, u)$ belongs to $\cT_X \Odr$.

Drawing on the optimal center $z\ustar_X \in \cC_X$ identified in Theorem~\ref{thm:center_attainment}, the following result provides an explicit characterization of the Riemannian Clarke subdifferential $\partial\lrmR \sfV_{n, \rho} (X)$. 

\begin{theorem}
	\label{thm:value_clarke_subdifferential}
	The function $\sfV_{n, \rho}$ is locally Lipschitz continuous. 
	Moreover, its Riemannian Clarke subdifferential is
	\begin{equation}
		\label{eq:value_clarke_subdifferential}
		\partial\lrmR \sfV_{n, \rho} (X)
		= \hkh{
			\sfG_X \dkh{X, z\ustar_X, u}
			\mid
			u \in \cU_0 \dkh{X, z\ustar_X},
			\sfG_z \dkh{X, z\ustar_X, u} = 0
		}.
	\end{equation}
\end{theorem}

\begin{proof}
	We fix $X \in \Odr$. 
	Let $\sfM (Y, z) := \sum_{i = 1}^n (\abs{Y_i - z}_2 + \rho)^2 / n$ for $Y = [Y_1, \dotsc, Y_n] \in \Rdn$ and $z \in \Rd$. 
	The function $\sfM$ is jointly convex in $(Y, z)$ and strongly convex in $z$. 
	Moreover, $\sfM$ is level-bounded in $z$ locally uniformly in $Y$, and its global minimizer $\tilde{z} (Y)$ is unique. 
	Then the marginal function $\sfW (Y) := \min_z \sfM (Y, z)$ is finite and convex on $\Rdn$. 
	It follows from Theorem~\ref{thm:center_attainment} that $\sfV_{n, \rho} (X) = \sfW (\sfY (X))$ and $\tilde{z} (\sfY (X)) = z\ustar_X$, where $\sfY(X) = [\sfY_1 (X), \dotsc, \sfY_n (X)]$ with $\sfY_i (X) = (I_d - X X\zz) \xi_i\ucirc$ for $i \in [n]$. 
	The convex subgradient inequality shows that $Q \in \partial\sfW (Y)|_{Y = \sfY(X)}$ if and only if $(Q, 0) \in \partial \sfM (Y, z)|_{(Y, z) = (\sfY(X), z\ustar_X)}$. 
	Then we can obtain that
	\begin{equation}
		\label{eq:subdifferential}
		\partial \sfW (Y) |_{Y = \sfY (X)}
		= \hkh{
			Q = [q_1, \dotsc, q_n] \;\middle|\;
			\begin{aligned}
				& q_i = 2 \dkh{\sfr_i (X, z\ustar_X) + \rho a_i} / n, \\
				& a_i \in \partial\norm{\cdot}_2 (x) |_{x = \sfr_i (X, z\ustar_X)},
				\sum\nolimits_{i = 1}^n q_i = 0
			\end{aligned}
		},
	\end{equation}
	where $\partial \norm{\cdot}_2 (x) = \{x / \norm{x}_2\}$ when $x \neq 0$, and the closed unit ball when $x = 0$. 
	Since $\sfW$ is finite and convex and $\sfY$ is smooth, the function $\sfV_{n,\rho} = \sfW \circ \sfY$ is locally Lipschitz continuous and Clarke regular. 
	The chain rule \cite[Theorem~10.6]{Rockafellar1998Variational} and the tangent-projection formula \cite[Theorem~5.1(iii)]{Yang2014Optimality} therefore yield
	\begin{equation*}
		\partial\lrmR \sfV_{n, \rho} (X)
		= \proj_{\cT_X \Odr} \dkh{\rmD \sfY(X)^\sharp \partial \sfW(Y) |_{Y = \sfY (X)}},
	\end{equation*}
	where $\rmD \sfY(X)^\sharp$ denotes the adjoint of the differential $\rmD \sfY (X)$. 
	
	Now we choose $Q \in \partial \sfW (Y) |_{Y = \sfY (X)}$ with the directions $a_i$ specified in \eqref{eq:subdifferential}. 
	Let $u_i := (I_d - X X\zz) a_i$. 
	Since $\sfr_i (X, z\ustar_X) \in \cS_X^\perp$, we have $u = (u_1, \ldots, u_n) \in \cU_0 (X, z\ustar_X)$ and $(I_d - X X\zz) q_i = 2 (\sfr_i (X, z\ustar_X) + \rho u_i) / n$.
	The condition $\sum_{i = 1}^n q_i = 0$ in \eqref{eq:subdifferential} thus implies that $\sfG_z (X, z\ustar_X, u) = 0$. 
	Furthermore, it is clear that $q_i \in \cS_X^\perp$ when $\sfr_i (X, z\ustar_X) \neq 0$ and $\sfY_i (X) = z\ustar_X$ when $\sfr_i (X, z\ustar_X) = 0$. 
	Then it can be readily verified that
	\begin{equation*}
		\sum_{i = 1}^n \sfY_i (X) q_i\zz X
		= z\ustar_X \dkh{\sum_{i = 1}^n q_i}\zz X
		= 0.
	\end{equation*}
	For any $H \in \cT_X \Odr$, the identities $\rmD \sfY_i (X) [H] = - (H X\zz + X H\zz) \xi_i\ucirc$ and $X\zz \rmD \sfY_i (X) [H] = - H\zz Y_i$ now give
	\begin{equation*}
		\begin{aligned}
			\jkh{\rmD \sfY(X)^\sharp Q, H}
			= {} & - \sum_{i = 1}^n \jkh{((I_d - X X\zz) q_i) (\xi_i\ucirc)\zz X, H}
			- \jkh{\sum_{i = 1}^n \sfY_i (X) q_i\zz X, H} \\
			= {} & \jkh{\sfG_X (X, z\ustar_X, u), H}.
		\end{aligned}
	\end{equation*}
	As $\sfG_X (X, z\ustar_X, u) \in \cT_X \Odr$, this proves the inclusion of $\partial\lrmR \sfV_{n, \rho} (X)$ in the right-hand side of \eqref{eq:value_clarke_subdifferential}. 
	Conversely, for any $u \in \cU_0 (X, z\ustar_X)$ satisfying $\sfG_z (X, z\ustar_X, u) = 0$, we choose $a_i = u_i$ and $q_i = 2 (\sfr_i (X, z\ustar_X) + \rho u_i) / n$.
	Then it holds that $Q \in \partial \sfW (\sfY (X))$. 
	The same calculation yields the reverse inclusion. 
	We complete the proof. 
\end{proof}


We say that $X \in \Odr$ is a Clarke stationary point if it satisfies $0 \in \partial\lrmR \sfV_{n, \rho} (X)$. 
Theorem~\ref{thm:value_clarke_subdifferential} establishes the variational foundation for analyzing the convergence of IRWPCA in terms of the original objective function $\sfV_{n, \rho}$. 
The subsequent results build on this characterization to develop a stationarity residual and derive theoretical guarantees.

\subsection{Stationarity Residual}

\label{subsec:value_stationarity_residual}

The exact characterization of $\partial\lrmR \sfV_{n, \rho}$ involves normalized residual directions, which may be discontinuous at zero sample residuals. 
For $\tau \geq 0$ and $X \in \Odr$, we introduce the following stationarity residual,
\begin{equation}
	\label{eq:value_stationarity_residual}
	\sfS_\tau^{\sfV} (X)
	:= \min_{z \in \cC_X, u \in \cU_\tau (X, z)}
	\dkh{
		\norm{\sfG_X (X, z, u)}\fs
		+
		\norm{\sfG_z (X, z, u)}_2^2
	}^{1 / 2},
\end{equation}
which relaxes the direction constraints at the prescribed scale $\tau$. 
The minimum of the above problem is attained because the feasible set is compact and the objective function is continuous. 

As shown below, the stationarity residual defined in \eqref{eq:value_stationarity_residual} vanishes precisely at Clarke stationary points when $\tau = 0$.

\begin{lemma}
	\label{lem:value_stationarity_residual_basic}
	A point $X \in \Odr$ is Clarke stationary if and only if $\sfS_0^{\sfV} (X) = 0$. 
\end{lemma}

\begin{proof}
	Suppose first that $\sfS_0^{\sfV} (X) = 0$.
	Then there exist $z \in \cC_X$ and $u \in \cU_0 (X, z)$ such that $\sfG_X (X, z, u) = 0$ and $\sfG_z (X, z, u) = 0$. 
	The equality $\sfG_z (X, z, u) = 0$ is the first-order condition for the convex minimization problem of $\sfF_{n, \rho} (X, z)$ with respect to $z \in \cS_X^\perp$ for fixed $X \in \Odr$. 
	The uniqueness of the optimal center in Theorem~\ref{thm:center_attainment} gives $z = z\ustar_X$. 
	It then follows from Theorem~\ref{thm:value_clarke_subdifferential} that $0 \in \partial\lrmR \sfV (X)$. 
	Conversely, we assume that $X \in \Odr$ is Clarke stationary. 
	Then there exist $u \in \cU_0 (X, z\ustar_X)$ such that both $\sfG_X (X, z\ustar_X, u)$ and $\sfG_z (X, z \ustar_X, u)$ vanish. 
	Taking $z = z\ustar_X$ in \eqref{eq:value_stationarity_residual} yields $\sfS_0^{\sfV} (X) = 0$. 
	The proof is completed.
\end{proof}

Since the feasible set in \eqref{eq:value_stationarity_residual} can only enlarge as the relaxation scale increases, we have $\sfS_{\tau_2}^{\sfV} (X) \leq \sfS_{\tau_1}^{\sfV} (X)$ if $0 \leq \tau_1 \leq \tau_2$. 
The following result further proves that this residual is consistent with Clarke stationarity as $\tau \to 0$.

\begin{lemma}
	\label{lem:value_residual_consistency}
	Suppose that there exist a sequence $\{X^{(j)}\} \subseteq \Odr$ and a nonnegative sequence $\{\tau^{(j)}\}$ such that $X^{(j)} \to \hat{X}$, $\tau^{(j)} \to 0$, and $\sfS_{\tau^{(j)}}^{\sfV} (X^{(j)}) \to 0$ as $j \to \infty$. 
	Then we have $0 \in \partial\lrmR \sfV (\hat{X})$.
\end{lemma}

\begin{proof}
	For each $j$, we choose $z^{(j)} \in \cC_{X^{(j)}}$ and $u_j \in \cU_{\tau^{(j)}} (X^{(j)}, z^{(j)})$ attaining the minimum in \eqref{eq:value_stationarity_residual}. 
	Then it holds that $( \sfS_{\tau^{(j)}}^{\sfV} (X^{(j)}) )^2 = \|\sfG_X (X^{(j)}, z^{(j)}, u^{(j)})\|\fs + \|\sfG_z (X^{(j)}, z^{(j)}, u^{(j)})\|_2^2$. 
	Since each $z^{(j)}$ is a convex combination of the projected samples $(I_d - X^{(j)} (X^{(j)})\zz) \xi_i\ucirc$, the sequence $\{z^{(j)}\}$ is bounded. 
	According to the definition in \eqref{eq:direction}, we also have $\|u_i^{(j)}\|_2 \leq 1$ for all $i$ and $j$.
	After passing to a common subsequence, we assume without loss of generality that $z^{(j)} \to \hat{z}$ and $u_i^{(j)} \to \hat{u}_i$ for $i \in [n]$ as $j \to \infty$. 
	The closedness of the simplex representation for $\cC_{X^{(j)}}$ gives $\hat{z} \in \cC_{\hat{X}}$. 
	Now we fix $i \in [n]$. If $\|\sfr_i (\hat{X}, \hat{z})\|_2 > 0$, it holds that $\|\sfr_i (X^{(j)}, z^{(j)})\|_2 > \tau_j$ for all sufficiently large $j$. 
	Then it can be readily verified that $\hat{u}_i = \sfr_i (\hat{X}, \hat{z}) / \|\sfr_i (\hat{X}, \hat{z})\|_2$. 
	If $\sfr_i (\hat{X}, \hat{z}) = 0$, we have $(X^{(j)})\zz u_i^{(j)} = 0$ and $\|u_i^{(j)}\|_2 \leq 1$ for each $j$. 
	Passing to the limit gives $\hat{u}_i \in \cS_{\hat{X}}^\perp$ and $\|\hat{u}_i\|_2 \leq 1$. 
	Hence, it holds that $\hat{u} \in \cU_0 (\hat{X}, \hat{z})$. 
	Both $\sfG_X$ and $\sfG_z$ in \eqref{eq:value_first_order_vectors} are continuous with respect to their arguments. 
	Since $\sfS_{\tau^{(j)}}^{\sfV} (X^{(j)}) \to 0$, we can obtain that $\sfG_X (X^{(j)}, z^{(j)}, u^{(j)}) \to 0$ and $\sfG_z (X^{(j)}, z^{(j)}, u^{(j)}) \to 0$, which further implies that $\sfG_X (\hat{X}, \hat{z}, \hat{u}) = 0$ and $\sfG_z (\hat{X}, \hat{z}, \hat{u}) = 0$. 
	It follows from $\hat{z} \in \cC_{\hat{X}}$ and $\hat{u} \in \cU_0 (\hat{X}, \hat{z})$ that $\sfS_0^{\sfV} (\hat{X}) = 0$. 
	By invoking Lemma~\ref{lem:value_stationarity_residual_basic}, we arrive at $0 \in \partial\lrmR \sfV (\hat{X})$, which completes the proof. 
\end{proof}

In view of the preceding two lemmas, the proposed stationarity residual provides a suitable measure for finite-accuracy analysis. 
Accordingly, for any $\epsilon > 0$, we say that a point $X \in \Odr$ is $\epsilon$-stationary if $\sfS_\epsilon^{\sfV} (X) \leq \epsilon$. 
This criterion requires a small stationarity residual while also restricting the relaxation scale.

\subsection{Stationarity Control via the Majorization Gap}

\label{subsec:convergence_analysis}

We next establish a quantitative link between the stationarity residual and the majorization gap introduced in \eqref{eq:majorization_gap}. 
This result connects the update mechanism of IRWPCA to the theoretical properties of problem~\eqref{opt:infwpca}.

The following lemma records a sufficient reduction of the merit function defined in \eqref{eq:merit}, with the reduction at each iteration bounded below by the corresponding majorization gap.

\begin{lemma}
	\label{lem:value_gap_summability}
	The iterates generated by Algorithm~\ref{alg:irwpca} satisfy
	\begin{equation}
		\label{eq:descent_merit}
		\sfE_{\tau^{(k)}} (X^{(k)}, z^{(k)}) 
		- \sfE_{\tau^{(k + 1)}} (X^{(k + 1)}, z^{(k + 1)})
		\geq \Delta^{(k)},
	\end{equation}
	for all $k \in \bN$. 
\end{lemma}

\begin{proof}
	According to Proposition~\ref{prop:joint_majorant_minimization} and \eqref{eq:merit_majorization_by_Q}, we have
	\begin{equation}
		\label{eq:sfqk}
		\sfQ^{(k)} (X^{(k + 1)}, z^{(k + 1)})
		\leq \sfQ^{(k)} (X^{(k)}, z^{(k)})
		= \sfE_{\tau^{(k)}} (X^{(k)}, z^{(k)}).
	\end{equation}
	Since $\tau^{(k + 1)} \leq \tau^{(k)}$, the definition in \eqref{eq:merit} yields $\sfE_{\tau^{(k + 1)}} (X, z) \leq \sfE_{\tau^{(k)}} (X, z)$ for all $(X, z) \in \Odr \times \Rd$. 
	Then it can be readily verified that
	\begin{equation*}
		\sfE_{\tau^{(k + 1)}} (X^{(k + 1)}, z^{(k + 1)}) + \Delta^{(k)}
		\leq \sfE_{\tau^{(k)}} (X^{(k + 1)}, z^{(k + 1)}) + \Delta^{(k)}
		= \sfQ^{(k)} (X^{(k + 1)}, z^{(k + 1)}),
	\end{equation*}
	which, together with \eqref{eq:sfqk}, implies the desired result \eqref{eq:descent_merit}. 
	The proof is completed.
\end{proof}

A key result is that the stationarity residual at the new iterate can be controlled by the normalized majorization gap, as established below.

\begin{proposition}
	\label{prop:stationarity_gap_bound}
	For the iterates generated by Algorithm~\ref{alg:irwpca}, we have
	\begin{equation}
		\label{eq:value_gap_residual_bound}
		\sfS_{\tau^{(k)}}^{\sfV} (X^{(k + 1)})
		\leq \sqrt{ \dfrac{4 \rho (\nu^2 + 1) \Delta^{(k)}}{\tau^{(k)}} },
	\end{equation}
	where $\nu = \max_{i \in [n]} \norm{\xi_i\ucirc}_2$ is a constant defined in Theorem~\ref{thm:center_attainment}. 
\end{proposition}

\begin{proof}
	For convenience, we define $\sfr_i^{(k + 1)} = \sfr_i (X^{(k + 1)}, z^{(k + 1)})$. 
	The update scheme of $z^{(k + 1)}$ in \eqref{eq:z_update} gives
	\begin{equation}
		\label{eq:optimality-z}
		\sum_{i = 1}^n \omega_i^{(k)} \sfr_i^{(k + 1)}
		= \dkh{I_d - X^{(k + 1)} (X^{(k + 1)})\zz} \sum_{i = 1}^n \omega_i^{(k)} (\xi_i\ucirc - z^{(k + 1)})
		= 0.
	\end{equation}
	By virtue of the optimality condition of problem \eqref{eq:X_update_eig}, we then arrive at
	\begin{equation}
		\label{eq:optimality-X}
		\sum_{i = 1}^n \omega_i^{(k)} \sfr_i^{(k + 1)} \dkh{\xi_i\ucirc}\zz X^{(k + 1)} = 0.
	\end{equation}
	Let $\tilde{\gamma}_i^{(k)} := \max \{\|\sfr_i^{(k + 1)}\|_2, \tau^{(k)}\}$ and $\tilde{u}_i^{(k)} = \sfr_i^{(k + 1)} / \tilde{\gamma}_i^{(k)}$ for $i \in [n]$. 
	Then it holds that
	\begin{equation*}
		\tilde{u}^{(k)} = [\tilde{u}_1^{(k)}, \dotsc, \tilde{u}_n^{(k)}] \in \cU_{\tau^{(k)}} (X^{(k + 1)}, \tilde{z}^{(k + 1)}),
	\end{equation*}
	where $\tilde{z}^{(k + 1)} := (I_d - X^{(k + 1)} (X^{(k + 1)})\zz) z^{(k + 1)} \in \cC_{X^{(k + 1)}}$.
	Since $\sfr_i^{(k + 1)} + \rho \tilde{u}_i^{(k)} = (1 + \rho / \tilde{\gamma}_i^{(k)}) \sfr_i^{(k + 1)}$ for all $i \in [n]$, we have
	\begin{equation*}
		\begin{aligned}
			\sfG_z (X^{(k + 1)}, \tilde{z}^{(k + 1)}, \tilde{u}^{(k)}) 
			= {} & - \dfrac{2}{n} \sum_{i = 1}^n \dkh{\sfr_i^{(k + 1)} + \rho \tilde{u}_i^{(k + 1)}}
			= - \dfrac{2}{n} \sum_{i = 1}^n \dkh{1 + \dfrac{\rho}{\tilde{\gamma}_i^{(k)}}} \sfr_i^{(k + 1)} \\
			= {} & \dfrac{2}{n} \sum_{i = 1}^n \omega_i^{(k)} \sfr_i^{(k + 1)} - \dfrac{2}{n} \sum_{i = 1}^n \dkh{1 + \dfrac{\rho}{\tilde{\gamma}_i^{(k)}}} \sfr_i^{(k + 1)} \\
			= {} & \dfrac{2 \rho}{n} \sum_{i = 1}^n \dkh{\dfrac{1}{\gamma_i^{(k)}} - \dfrac{1}{\tilde{\gamma}_i^{(k)}}} \sfr_i^{(k + 1)},
		\end{aligned}
	\end{equation*}
	where the third equality uses \eqref{eq:optimality-z}. 
	Similarly, it follows from \eqref{eq:optimality-X} that
	\begin{equation*}
		\sfG_X (X^{(k + 1)}, \tilde{z}^{(k + 1)}, \tilde{u}^{(k)})
		= \dfrac{2 \rho}{n} \sum_{i = 1}^n \dkh{ \dfrac{1}{\gamma_i^{(k)}} - \dfrac{1}{\tilde{\gamma}_i^{(k)}} } \sfr_i^{(k + 1)} \dkh{\xi_i\ucirc}\zz X^{(k + 1)}.
	\end{equation*}
	Now we define
	\begin{equation*}
		\Gamma^{(k)} := \dfrac{1}{n} \sum_{i = 1}^n \dkh{ \dfrac{1}{\gamma_i^{(k)}} - \dfrac{1}{\tilde{\gamma}_i^{(k)}} }^2 \norm{\sfr_i^{(k + 1)}}_2^2.
	\end{equation*}
	Applying the triangle and Cauchy-Schwarz inequalities yields that
	\begin{equation*}
		\norm{ \sfG_z (X^{(k + 1)}, \tilde{z}^{(k + 1)}, \tilde{u}^{(k)}) }_2
		\leq \dfrac{2 \rho}{n} \sum_{i = 1}^n \abs{ \dfrac{1}{\gamma_i^{(k)}} - \dfrac{1}{\tilde{\gamma}_i^{(k)}} } \norm{ \sfr_i^{(k + 1)} }_2
		\leq 2 \rho \sqrt{ \Gamma^{(k)} }.
	\end{equation*}
	By similar arguments, we can obtain that
	\begin{equation*}
		\norm{ \sfG_X (X^{(k + 1)}, \tilde{z}^{(k + 1)}, \tilde{u}^{(k)}) }\ff
		\leq 2 \rho \nu \sqrt{ \Gamma^{(k)} },
	\end{equation*}
	which further implies that
	\begin{equation}
		\label{eq:residual}
		\sfS_{\tau^{(k)}}^{\sfV} (X^{(k + 1)})
		\leq 2 \rho \sqrt{(\nu^2 + 1) \Gamma^{(k)}}. 
	\end{equation}
	It remains to prove that $\Gamma^{(k)} \leq \Delta^{(k)} / (\rho \tau^{(k)})$. 
	Let
	\begin{equation*}
		\theta_i^{(k)}
		= \dfrac{1}{\gamma_i^{(k)}} \norm{\sfr_i^{(k + 1)}}_2^2 
		+ \gamma_i^{(k)} 
		- \dfrac{1}{\tilde{\gamma}_i^{(k)}} \norm{\sfr_i^{(k + 1)}}_2^2
		- \tilde{\gamma}_i^{(k)},
	\end{equation*}
	for $i \in [n]$. 
	The definitions of $\sfQ^{(k)}$ and $\sfE_{\tau^{(k)}}$ give $\Delta^{(k)} = \rho \sum_{i = 1}^{n} \theta_i^{(k)} / n$. 
	If $\|\sfr_i^{(k + 1)}\|_2 \geq \tau^{(k)}$, we have $\tilde{\gamma}_i^{(k)} := \|\sfr_i^{(k + 1)}\|_2$ and $\theta_i^{(k)} = (\gamma_i^{(k)} - \|\sfr_i^{(k + 1)}\|_2)^2 / \gamma_i^{(k)}$. 
	If $\|\sfr_i^{(k + 1)}\|_2 < \tau^{(k)}$, it holds that $\tilde{\gamma}_i^{(k)} := \tau^{(k)}$ and
	\begin{equation*}
		\theta_i^{(k)}
		= \dkh{ \gamma_i^{(k)} - \tau^{(k)} } \dkh{ 1 - \dfrac{1}{\gamma_i^{(k)} \tau^{(k)}} \norm{\sfr_i^{(k + 1)}}_2^2 }
		\geq \dfrac{\dkh{\gamma_i^{(k)} - \tau^{(k)}}^2}{\gamma_i^{(k)}}.
	\end{equation*}
	In both cases, straightforward calculations reveal that
	\begin{equation*}
		\dkh{ \dfrac{1}{\gamma_i^{(k)}} - \dfrac{1}{\tilde{\gamma}_i^{(k)}} }^2 \norm{\sfr_i^{(k + 1)}}_2^2
		\leq \dfrac{\theta_i^{(k)}}{\tau^{(k)}}.
	\end{equation*}
	The inequality $\Gamma^{(k)} \leq \Delta^{(k)} / (\rho \tau^{(k)})$ then follows by summing the above relationship over $i$. 
	As a direct consequence of \eqref{eq:residual}, we finally obtain the estimate \eqref{eq:value_gap_residual_bound}. 
	The proof is completed. 
\end{proof}

Combining the previous results with the selection rule of IRWPCA, we obtain an explicit bound on the stationarity residual at the selected iterates.

\begin{corollary}
	\label{cor:sublinear_stationarity_general}
	The output iterates selected by Algorithm~\ref{alg:irwpca} satisfy
	\begin{equation}
		\label{eq:general_complexity_bound}
		\sfS_{\bar{\tau}^{(k)}}^{\sfV} (\bar{X}^{(k + 1)})
		\leq \sqrt{ \dfrac{C}{(k + 1)^{1 - \sigma}} },
	\end{equation}
	where $C := 8 \rho (\nu^2 + 1) \sfE_{\tau^{(0)}} (X^{(0)}, z^{(0)}) / \tau^{(0)}$ is a constant.
\end{corollary}

\begin{proof}
	Let $\cJ^{(k)} := \{\floor{(k + 1) / 2}, \dotsc, k\}$. 
	This index set contains at least $(k + 1) / 2$ elements. 
	Summing the relationship \eqref{eq:descent_merit} yields
	\begin{equation*}
		\begin{aligned}
			\sum_{j \in \cJ^{(k)}} \Delta^{(j)}
			\leq {} & \sum_{j = 0}^{k} \dkh{ \sfE_{\tau^{(j)}} (X^{(j)}, z^{(j)}) - \sfE_{\tau^{(j + 1)}} (X^{(j + 1)}, z^{(j + 1)}) } \\
			= {} & \sfE_{\tau^{(0)}} (X^{(0)}, z^{(0)}) - \sfE_{\tau^{(k + 1)}} (X^{(k + 1)}, z^{(k + 1)}) \\
			\leq {} & \sfE_{\tau^{(0)}} (X^{(0)}, z^{(0)}),
		\end{aligned}
	\end{equation*}
	where the last inequality uses $\sfE_{\tau^{(k + 1)}} (X^{(k + 1)}, z^{(k + 1)}) \geq 0$. 
	The safeguard schedule guarantees that $\tau^{(j)} \geq \tau^{(k)}$ for every $j \in \cJ^{(k)}$. 
	According to the definition of $\sfm (k)$ in Algorithm~\ref{alg:irwpca}, it follows that
	\begin{equation*}
		\dfrac{\Delta^{(\sfm (k))}}{\tau^{(\sfm (k))}}
		= \min_{j \in \cJ^{(k)}} \dfrac{\Delta^{(j)}}{\tau^{(j)}}
		\leq \dfrac{2}{k + 1} \sum_{j \in \cJ^{(k)}} \dfrac{\Delta^{(j)}}{\tau^{(j)}}
		\leq \dfrac{2}{(k + 1) \tau^{(k)}} \sum_{j \in \cJ^{(k)}} \Delta^{(j)}
		\leq \dfrac{2 \sfE_{\tau^{(0)}} (X^{(0)}, z^{(0)})}{(k + 1) \tau^{(k)}}.
	\end{equation*}
	By Proposition~\ref{prop:stationarity_gap_bound}, we proceed to show that
	\begin{equation*}
		\sfS_{\bar{\tau}^{(k)}}^{\sfV} (\bar{X}^{(k + 1)})
		= \sfS_{\tau^{(\sfm (k))}}^{\sfV} (X^{(\sfm (k) + 1)})
		\leq \sqrt{ \dfrac{4 \rho (\nu^2 + 1) \Delta^{(\sfm (k))}}{\tau^{(\sfm (k))}} }
		\leq
		\sqrt{ \dfrac{ C \tau^{(0)} }{ (k + 1) \tau^{(k)} } }.
	\end{equation*}
	Substituting $\tau^{(k)} = \tau^{(0)} (k + 1)^{- \sigma}$ into the above inequality completes the proof. 
\end{proof}

With these preparations in place, we are in a position to develop a comprehensive convergence theory of our algorithm.

\subsection{Global Convergence and Iteration Complexity}

\label{subsec:epsilon_complexity}

The variational characterization and majorization-gap estimates now culminate in the convergence guarantees for IRWPCA. 
From any initialization, every accumulation point of the generated sequence is shown to be Clarke stationary. 
We also derive an iteration-complexity bound for obtaining an $\epsilon$-stationary point.

\begin{theorem}
	\label{thm:epsilon_complexity}
	Let $\{\bar{X}^{(k)}\}$ be the output sequence generated by Algorithm~\ref{alg:irwpca}. 
	Then this sequence admits at least one accumulation point, every one of which is Clarke stationary. 
	Furthermore, the selected output $\bar{X}^{(k_\epsilon)}$ qualifies as an $\epsilon$-stationary point, where
	\begin{equation*}
		k_\epsilon
		= \ceil{
			\max \hkh{
				1,
				\dfrac{ C^{1 / (1 - \sigma)} }{ \epsilon^{2 / (1 - \sigma)} },
				\dfrac{ 2 (\tau^{(0)})^{1 / \sigma} }{ \epsilon^{1 / \sigma} }
			}
		},
	\end{equation*}
	for any $\epsilon > 0$. 
\end{theorem}

\begin{proof}
	The compactness of the Stiefel manifold guarantees the existence of accumulation points. 
	Since $\bar{\tau}^{(k)} = \tau^{(\sfm (k))} \leq \tau^{(\floor{(k + 1) / 2})} = \tau^{(0)} (\floor{(k + 1) / 2} + 1)^{- \sigma} \leq 2^{\sigma} \tau^{(0)} (k + 1)^{- \sigma}$, we have $\bar{\tau}^{(k)} \to 0$ as $k \to \infty$. 
	Moreover, it follows from corollary~\ref{cor:sublinear_stationarity_general} that $\sfS_{\bar{\tau}^{(k)}}^{\sfV} (\bar{X}^{(k + 1)}) \to 0$ as $k \to \infty$. 
	Then Lemma~\ref{lem:value_residual_consistency} applies to any convergent subsequence of $\{\bar{X}^{(k)}\}$ and yields the Clarke stationarity of every accumulation point. 
	Combining the definition of $k_\epsilon$ with \eqref{eq:general_complexity_bound}, we then arrive at
	\begin{equation*}
		\sfS_{\bar{\tau}^{(k_\epsilon - 1)}}^{\sfV} (\bar{X}^{(k_\epsilon)})
		\leq \sqrt{ \dfrac{C}{k_\epsilon^{1 - \sigma}} }
		\leq \epsilon,
	\end{equation*}
	which, together with $\bar{\tau}^{(k_\epsilon - 1)} \leq 2^{\sigma} \tau^{(0)} k_\epsilon^{- \sigma} \leq \epsilon$, implies that $\sfS_\epsilon^{\sfV} (\bar{X}^{(k_\epsilon)}) \leq \sfS_{\bar{\tau}^{(k_\epsilon - 1)}}^{\sfV} (\bar{X}^{(k_\epsilon)}) \leq \epsilon$. 
	Therefore, the selected output $\bar{X}^{(k_\epsilon)}$ is an $\epsilon$-stationary point. 
	The proof is completed.
\end{proof}

Theorem~\ref{thm:epsilon_complexity} ensures that IRWPCA attains an $\epsilon$-stationary point within $O (\epsilon^{- \max \{2 / (1 - \sigma), 1 / \sigma\}})$ iterations for $\sigma \in (0, 1)$. 
In particular, the choice $\sigma = 1 / 3$ balances the two terms in the complexity exponent and leads to an $O (\epsilon^{- 3})$ iteration-complexity bound.

\section{Numerical Study}

\label{sec:numerical}

Preliminary numerical results are presented in this section to investigate the computational efficient of the IRWPCA algorithm and the practical effectiveness of the InfWPCA model in downstream applications. 
The corresponding experiments are conducted using Python on a MacBook Pro equipped with an Apple M2 Pro chip and $16$ GB of memory.

\subsection{Computational Behavior of IRWPCA}

\label{subsec:numerical_algorithmic}

The first experiment is designed to evaluate the numerical performance of IRWPCA in solving problem~\eqref{opt:infwpca-eq}, an exact reformulation of the InfWPCA model. 
We compare IRWPCA against the Riemannian subgradient method (RSG) of \cite{Li2021Weakly} and the dynamic smoothing gradient method (DSGM) of \cite{Beck2023Dynamic}. 
The two benchmark algorithms are adapted to handle the joint variables $(X, z)$ of problem~\eqref{opt:infwpca-eq} by treating $\Odr \times \Rd$ as a product manifold.

We randomly select a balanced sample set of $n = 3000$ images in the MNIST dataset \cite{LeCun1998Gradient}, with $300$ images from each of the ten classes. 
Every image contains $d = 784$ pixels, and the pixel values are normalized to $[0, 1]$ through division by $255$. 
In problem~\eqref{opt:infwpca-eq}, we set $\rho = 1$ and consider $r \in \{10, 50, 100\}$. 
All algorithms are initialized at the same point, which is randomly generated for each trial. 
For our testing, IRWPCA is configured with $\tau^{(0)} = 1 / 2$ and $\sigma = 1 / 3$ by default. 
In the implementation of RSG, the stepsizes follow $\eta\lRSG^{(k)} = \eta\lRSG^{(0)} (k + 1)^{- 1 / 2}$, with $\eta\lRSG^{(0)} \in \{0.003, 0.01, 0.03, 0.1, 0.3\}$. 
We employ a polar retraction and a deterministic Clarke subgradient selection. 
The smoothing parameters of DSGM are set to $\mu\lDSGM^{(k)} = 0.5 (k + 1)^{- 1 / 3}$, with initial stepsizes chosen from $\{0.03, 0.1, 0.3, 1.0\}$. 
Its Armijo and backtracking constants are fixed at $10^{- 4}$ and $0.5$, respectively. 
For the two benchmark algorithms, we retain the stepsize candidate that delivers the best performance in each test instance. 
By contrast, IRWPCA does not require a grid search for stepsizes.

\begin{figure}[ht]
	\centering
	
	\subfloat[$r = 10$]{
		\label{subfig:mnist_r10}
		\includegraphics[width=0.31\linewidth]{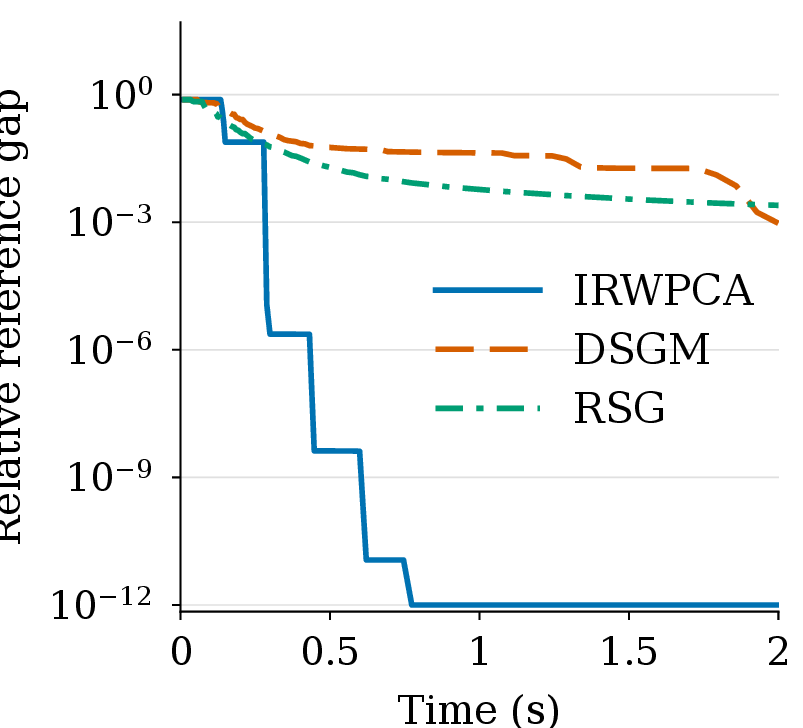}
	}
	\subfloat[$r = 50$]{
		\label{subfig:mnist_r50}
		\includegraphics[width=0.31\linewidth]{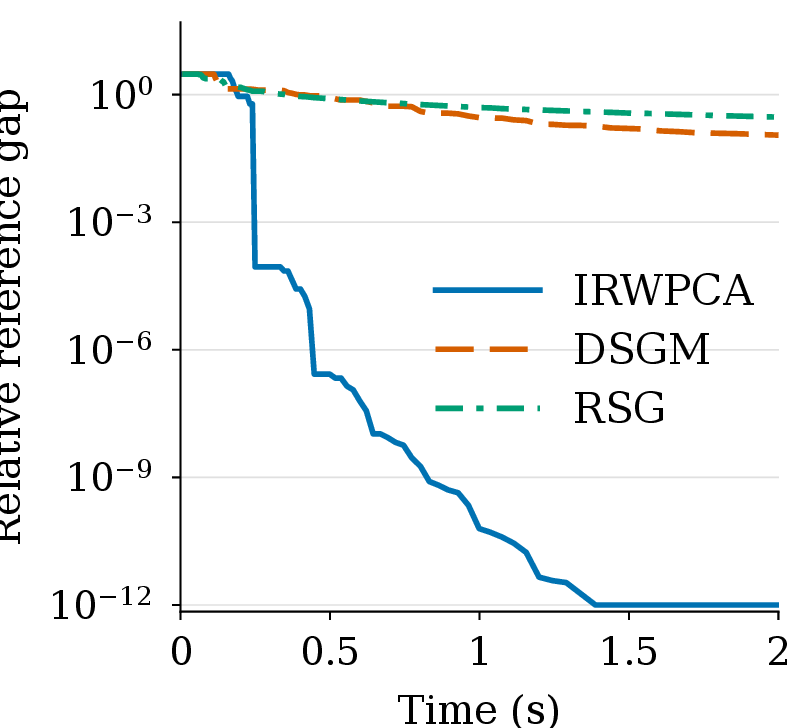}
	}
	\subfloat[$r = 100$]{
		\label{subfig:mnist_r100}
		\includegraphics[width=0.31\linewidth]{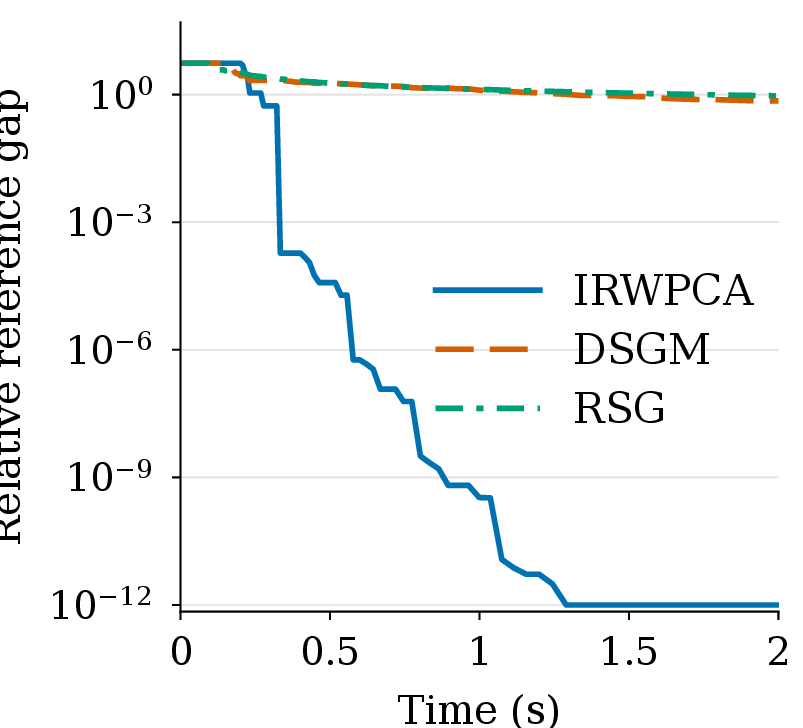}
	}
	
	\caption{Comparison of IRWPCA, DSGM, and RSG for solving the InfWPCA model.}
	\label{fig:numerical_algorithmic_cpu}
\end{figure}

To ensure a fair comparison, we impose a uniform time budget of $2$ seconds on all algorithms. 
The computational time is measured in wall-clock seconds, excluding data preparation and initialization. 
For each test instance, we conduct ten independent trials. 
The following relative reference gap is adopted as the performance metric,
\begin{equation*}
	\dfrac{\sfF_{n, \rho} (X, z) - \sfF\lbest}{\max \hkh{\abs{\sfF\lbest}, 10^{- 12}}} + 10^{- 12},
\end{equation*}
where $\sfF\lbest$ records the smallest objective value achieved by any of the tested algorithms within the prescribed budget in each trial. 
Figure~\ref{fig:numerical_algorithmic_cpu} illustrates the decrease in the relative reference gap against wall-clock time, with each curve averaged over ten trials. 
As evidenced therein, IRWPCA stands out in computational efficiency, attaining the lowest final objective value among the three algorithms across all cases. 
Moreover, the performance advantage of IRWPCA becomes more substantial as the value of $r$ increases. 
These observations support the numerical superiority of IRWPCA.

\subsection{Image Representation and Classification}

\label{subsec:image}

In this subsection, we assess the ability of InfWPCA to learn reliable low-dimensional representations in the presence of data contamination. 
The corresponding experiments are conducted on the $8 \times 8$ handwritten-digit dataset in scikit-learn \cite{Pedregosa2011Scikit}, where the pixel values are normalized to $[0, 1]$. 
We randomly generate ten stratified splits indexed, each consisting of $n = 1168$ training samples and $m = 629$ test samples. 
The contamination is introduced only into the training samples.

Let $\{\xi_i\utrain\}_{i = 1}^{n}$ and $\{\xi_i\utest\}_{i = 1}^{m}$ be the clean training and test samples, respectively. 
For each instance, a uniformly selected set of $\floor{p n}$ training samples is corrupted, where $p \in \{0.15, 0.30, 0.50\}$ is the contamination fraction. 
These fraction values correspond to $175$, $350$, and $584$ corrupted samples. 
The notation $\xi_i\ubad$ represents the $i$-th training sample after the contamination procedure, with unselected samples remaining unchanged. 
We apply the following three contamination mechanisms with magnitude $m > 0$ to the selected training samples. 
\begin{equation*}
	\left\{
	\begin{aligned}
		& \mbox{Gross noise:}
		&& \xi_i\ubad = \xi_i\utrain + m \eta_i,
		\quad \eta_i \sim \dN (0, I_d). \\
		& \mbox{Replacement:} 
		&& \xi_i\ubad = \mu\ltrain + m \phi\ltrain \odot \eta_i,
		\quad \eta_i \sim \dN (0, I_d). \\
		&\mbox{Sparse impulse:} 
		&& \xi_i\ubad = \xi_i\utrain + m u_i,
		\quad \norm{u_i}_0 = \floor{d / 4},
		\quad u_{i, j} \in \hkh{\pm1}.
	\end{aligned}
	\right.
\end{equation*}
Here, $\mu\ltrain = \sum_{i = 1}^{n} \xi_i\utrain / n$ is the mean vector of clean training samples, $\phi\ltrain \in \Rd$ contains the coordinatewise standard deviations of the clean training data, $\odot$ denotes the Hadamard product, $\dN$ represents the normal distribution, $\norm{u_i}_0$ returns the number of nonzero entries in $u_i$, and $u_{i, j}$ stands for the $j$-th entry of $u_i$. 
Unless stated otherwise, we fix $m = 3.5$ in the following experiments.

All the tested models are evaluated with $r = 10$ on the same data for a given split and contamination pattern. 
We deploy the IRWPCA algorithm for $100$ iterations to solve the InfWPCA model. 
The radius $\rho$ of the ambiguity set is selected by a validation strategy. 
Specifically, the contaminated training data is partitioned into a fitting set of $75\%$ samples and a validation set of $25\%$ samples. 
For each candidate multiplier $c \in \{10^{- 4 + 8 j / 19} \mid j = 0, 1, \dotsc, 19\}$, we solve problem~\eqref{opt:infwpca-eq} on the fitting set with $\rho = c / \sqrt{n}$. 
The median of reconstruction errors on the validation samples is evaluated as the performance metric. 
We then select the best-performing multiplier and refit the InfWPCA model on the full training set. 

In addition to classical PCA, we include two data-adaptive order-$2$ Wasserstein surrogates proposed in \cite{Xu2026Distributionally} for comparison. 
The corresponding model is referred to as W2PCA. 
Both surrogates minimizes the objective function in \eqref{eq:w2_upper_bound} over the Stiefel manifold, differing only in the choice of the weight matrix $G$. 
Let $A\lbad = [\xi_1\ubad - \mu\lbad, \dotsc, \xi_n\ubad - \mu\lbad] \in \Rdn$ with $\mu\lbad = \sum_{i = 1}^{n} \xi_i\ubad / n$ and $X\lpca \in \Odr$ be the solution of classical PCA. 
We also denote $\Sigma\lbad = A\lbad A\lbad\zz / n$ and $Q\lpca = X\lpca X\lpca\zz$. 
The residual-based and PCA-based variance vectors are defined by $v\lres = \diag ( (I_d - Q\lpca) \Sigma\lbad (I_d - Q\lpca) )$ and $v\lpca = \diag (Q\lpca \Sigma\lbad Q\lpca)$, respectively. 
A ridge term of $0.01 \tr (\Sigma\lbad) / d$ is added to each entry of the two variance vectors. 
Then we normalize each resulting vector to have mean one, which yields $b\lres$ and $b\lpca$. 
The residual-based and PCA-based surrogates employ the weight matrices $G\lres = \Diag (b\lres)\inv$ and $G\lpca = \Diag (b\lpca)\inv$, respectively. 
Xu et al.~\cite{Xu2026Distributionally} follow the robust Wasserstein profile inference framework of \cite{Blanchet2019Robust} to calibrate the radius. 
In our experiments, we set $\rho = \sqrt{q_{0.9}\ucirc / n}$ for each surrogate of W2PCA. 
Here, $q_{0.9}\ucirc$ is the estimated $0.9$ quantile of the limiting weighted chi-squared distribution, obtained from $4096$ Gaussian draws. 
We use the path-based grid implementation proposed in \cite{Xu2026Distributionally} to approximately solve the two surrogate problems. 
For each weight matrix $G$, candidate solutions are constructed from the $r$ leading eigenvectors of $\Sigma\lbad + \iota \rho \tr (\Sigma\lbad) G\inv / d$, where $\iota$ ranges over $61$ logarithmically spaced values in $[10^{- 3}, 10^3]$. 
The candidate set also includes the solution of classical PCA and the matrix formed by the $r$ standard basis vectors associated with the largest diagonal entries of $G\inv$. 
The candidate with the smallest value of the corresponding surrogate is then selected.

To evaluate the quality of a learned pair $(X, z)$, we calculate the mean squared error (MSE) of the reconstruction residuals on clean test samples as follows,
\begin{equation*}
	\MSE\ltest (X, z) = \dfrac{1}{m} \sum_{i = 1}^{m} \norm{(I_d - X X\zz) (\xi_i\utest - z)}_2^2.
\end{equation*}
The discriminative information preserved in the learned representations is further assessed through a classification task. 
Specifically, we train an $\ell_2$-regularized multinomial logistic regression classifier on the low-dimensional representations of contaminated training samples and evaluate its accuracy on the corresponding representations of clean test samples. 
All classifiers use a common regularization parameter of $1$ and the same LBFGS solver with a maximum of $1000$ iterations.

\begin{figure}[ht]
	\centering
	
	\subfloat[MSE]{
		\label{subfig:digits_mse}
		\includegraphics[width=0.48\linewidth]{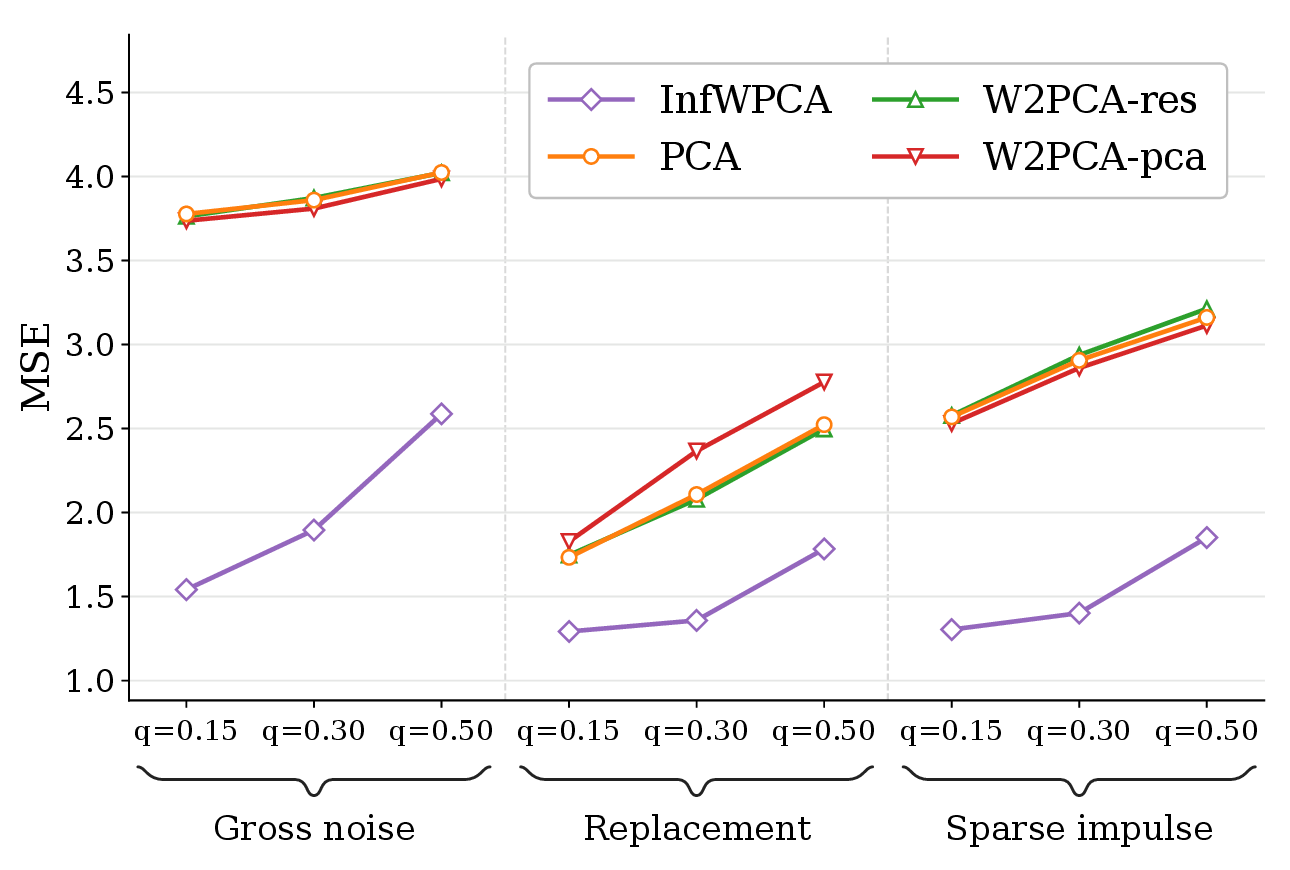}
	}
	\subfloat[Accuracy]{
		\label{subfig:digits_acc}
		\includegraphics[width=0.48\linewidth]{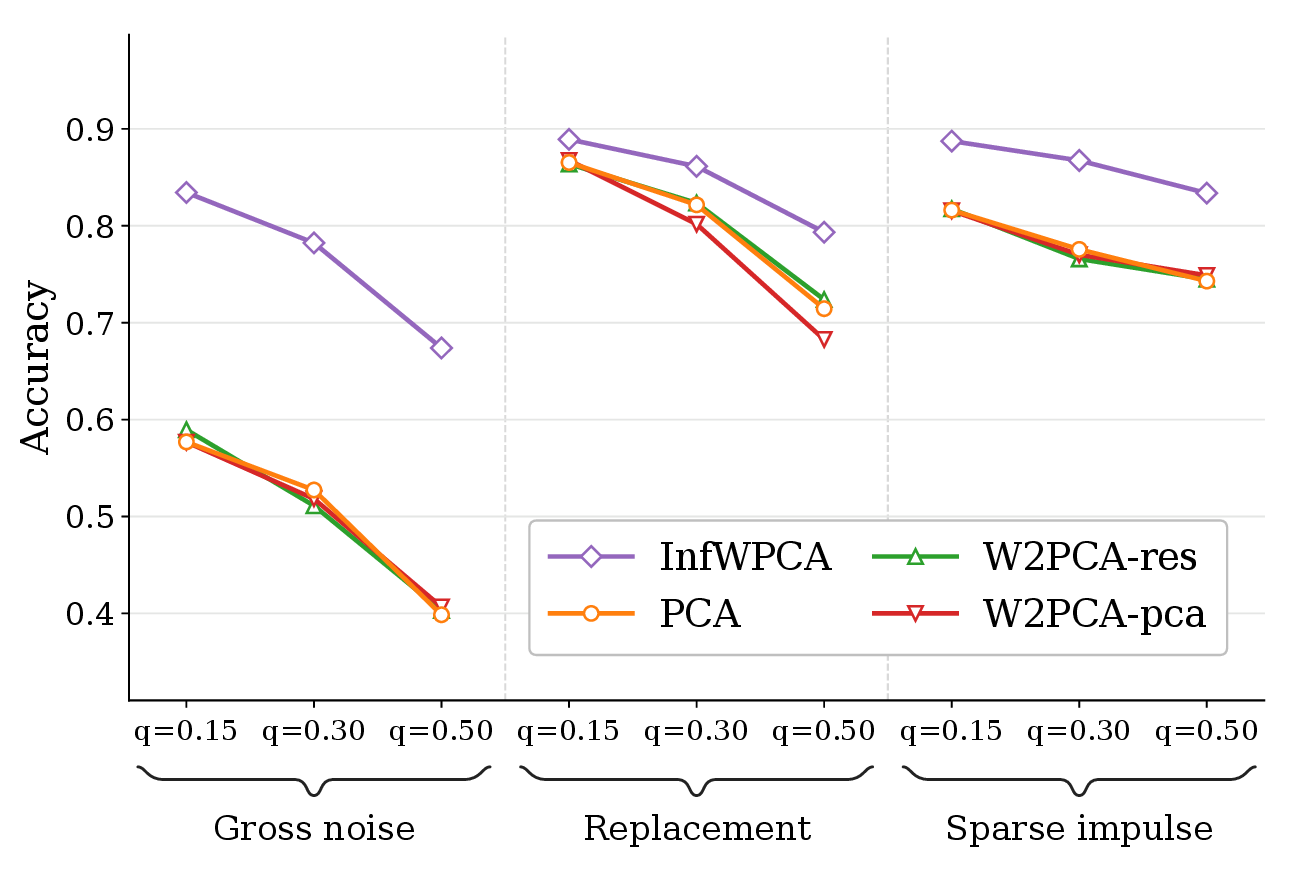}
	}
	
	\caption{Reconstruction MSE and classification accuracy on clean test samples. Smaller MSE and larger accuracy indicate better performance.}
	\label{fig:numerical_digits}
\end{figure}

Figure~\ref{fig:numerical_digits} presents the average reconstruction MSE and classification accuracy on clean test samples over ten data splits. 
W2PCA-res and W2PCA-pca refer to the W2PCA models with residual-based and PCA-based surrogates, respectively. 
Across all tested settings, InfWPCA delivers a significant improvement in both criteria, achieving the smallest reconstruction error and the highest classification accuracy. 
We can also observe that this performance gap is particularly pronounced under the Gaussian noise contamination. 
These numerical results demonstrate the favorable empirical performance of InfWPCA in learning low-dimensional representations from contaminated data.

\subsection{Minimum-Variance Portfolio Selection}

\label{subsec:numerical_finance}

Finally, we explore the practical impact of the subspace learned by InfWPCA on the downstream task of portfolio selection. 
These experiments are inspired by PCA-based factor models for covariance estimation and portfolio optimization \cite{Fan2008High,Fan2013Large,Kao2014Directed}. 
We use the monthly returns in the 48 Industry Portfolios dataset from the Fama-French library\footnote{\url{https://mba.tuck.dartmouth.edu/pages/faculty/ken.french/data_library.html}}. 
The input panel comprises $674$ consecutive monthly observations, from November 1969 through December 2025, across $d = 48$ industries. 
The initial $50$ months form the first estimation window, and the subsequent $624$ months, spanning January 1974 to December 2025, constitute the common out-of-sample evaluation period.

At each potential decision date, we fit the covariance estimator based on the $n = 50$ immediately preceding monthly return vectors. 
Let $l \in \{1, \dotsc, 12\}$ denote the portfolio update interval, measured in months. 
Starting in January 1974, we determine an allocation from the available $50$-month history, which ia applied to the ensuing $l$ monthly return vectors and replaced by a newly estimated allocation at the next scheduled update. 
The trailing estimation window is advanced by $l$ months at every subsequent re-estimation. 
All twelve choices of $l$ are evaluated over the same $T = 624$ test months. 
Only the frequency with which the portfolio is recalibrated changes.

We consider the four models examined in Section~\ref{subsec:image} with $r = 5$, including classical PCA, InfWPCA, W2PCA-res and W2PCA-pca. 
The same parameter-tuning protocols and implementation settings are adopted here, except that all models are fitted to the original data without introducing contamination. 
Given a learned pair $(X, z)$, we reconstruct the training samples as
\begin{equation*}
	R\lrec = X X\zz R\ltrain + \dkh{ I_d - X X\zz} z \bfone_n\zz,
\end{equation*}
where $\bfone_n$ is the $n$-dimensional vector of all ones and $R\ltrain \in \Rdn$ contains the current training returns. 
Then the reconstructed covariance matrix is given by
\begin{equation*}
	\Sigma\lrec = \dfrac{1}{n - 1} \dkh{R\lrec - R\lrec \bfone_n \bfone_n\zz / n} \dkh{R\lrec - R\lrec \bfone_n \bfone_n\zz / n}\zz.
\end{equation*}
To improve numerical stability, we apply a common ridge adjustment to the reconstructed covariance matrix, yielding $\Sigma\lport = \Sigma\lrec + 10^{- 4} \tr (\Sigma\lrec) I_d / d$. 
The estimated covariance matrix is subsequently employed in the long-only minimum-variance portfolio optimization problem \cite{Markowitz1952Portfolio},
\begin{equation}
	\label{opt:portfolio}
	\min_{w \geq 0, \bfone_d\zz w = 1} w\zz \Sigma\lport w.
\end{equation}
We use the same SLSQP solver to tackle this problem with a maximum of $900$ iterations, a termination tolerance of $10^{- 12}$, and explicit feasibility checks.

For a given update interval $l$, the estimated weight vector of problem~\eqref{opt:portfolio} is kept constant between successive update dates and applied to the realized industry returns in each intervening month. 
The backtest thereby produces $T$ monthly portfolio returns $\{a_1, \dotsc, a_T\}$. 
Let $\bar{a}_T =\sum_{t = 1}^T a_t / T$ and $s_T^2 = \sum_{t = 1}^T (a_t - \bar{a}_T)^2 / (T - 1)$ be the mean and variance for one such return series. 
We assess the out-of-sample performance of the tested models based on four metrics, including annualized mean return (AMR), Sharpe ratio (SR), annualized volatility (AV), and maximum drawdown (MD). 
These metrics are computed as follows,
\begin{equation*}
	\AMR = 12 \bar{a}_T,
	\quad
	\SR = \dfrac{\sqrt{12} \bar{a}_T}{s_T},
	\quad
	\AV = \sqrt{12} s_T,
	\quad
	\MD = \max_{1\leq t \leq T} \dkh{ 1 - \dfrac{W_t}{\max_{0 \leq u \leq t} W_u} },
\end{equation*}
where $W_0 = 1$ and $W_t = \prod_{u = 1}^{t} (1 + a_u)$ for $t \geq 1$. 
The monthly risk-free rate is set to zero in the calculation of Sharpe ratios. 
Higher returns and Sharpe ratios are desirable, whereas lower volatility and drawdown signify greater risk control.

\begin{table}[ht]
	\centering
	\caption{Averaged values of the four performance metrics for the tested models over twelve update frequencies. The best results are highlighted in boldface.}
	\label{tab:portfolio}
	\small
	\setlength{\tabcolsep}{5pt}
	\begin{tabular}{lcccc}
		\toprule
		Model & SR $\uparrow$ & AMR $\uparrow$ & AV $\downarrow$ & MD $\downarrow$ \\
		\midrule
		InfWPCA & \textbf{0.835} & \textbf{11.435\%} & \textbf{13.700\%} & \textbf{38.418\%} \\
		PCA & 0.816 & 11.336\% & 13.886\% & 41.811\% \\
		W2PCA-res & 0.760 & 10.473\% & 13.786\% & 40.447\% \\
		W2PCA-pca & 0.783 & 11.145\% & 14.233\% & 38.642\% \\
		\bottomrule
	\end{tabular}
\end{table}

\begin{figure}[ht]
	\centering
	
	\subfloat[Sharpe Ratio]{
		\label{subfig:port_sr}
		\includegraphics[width=0.48\linewidth]{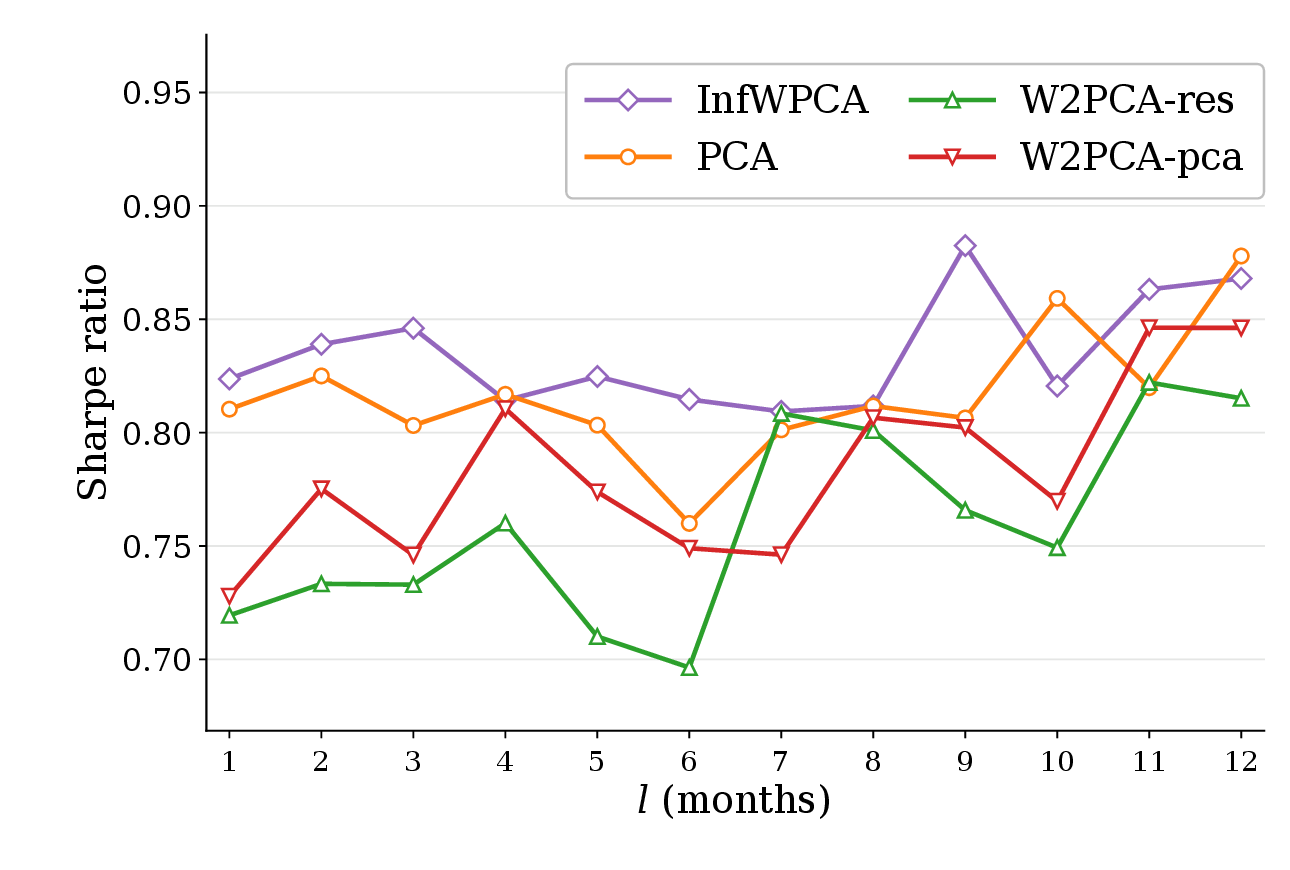}
	}
	\subfloat[Annualized Mean Return]{
		\label{subfig:port_amr}
		\includegraphics[width=0.48\linewidth]{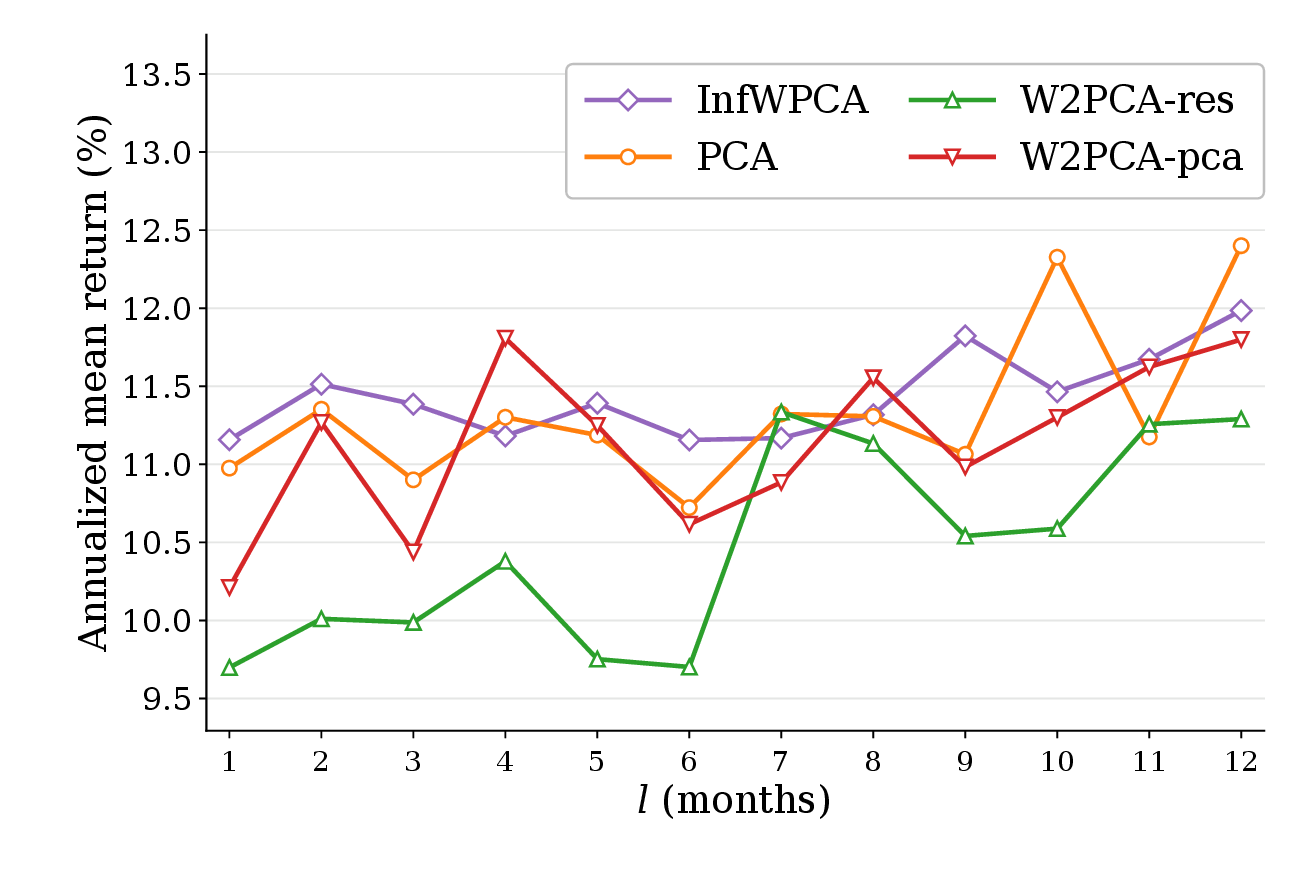}
	}
	
	\subfloat[Annualized Volatility]{
		\label{subfig:port_av}
		\includegraphics[width=0.48\linewidth]{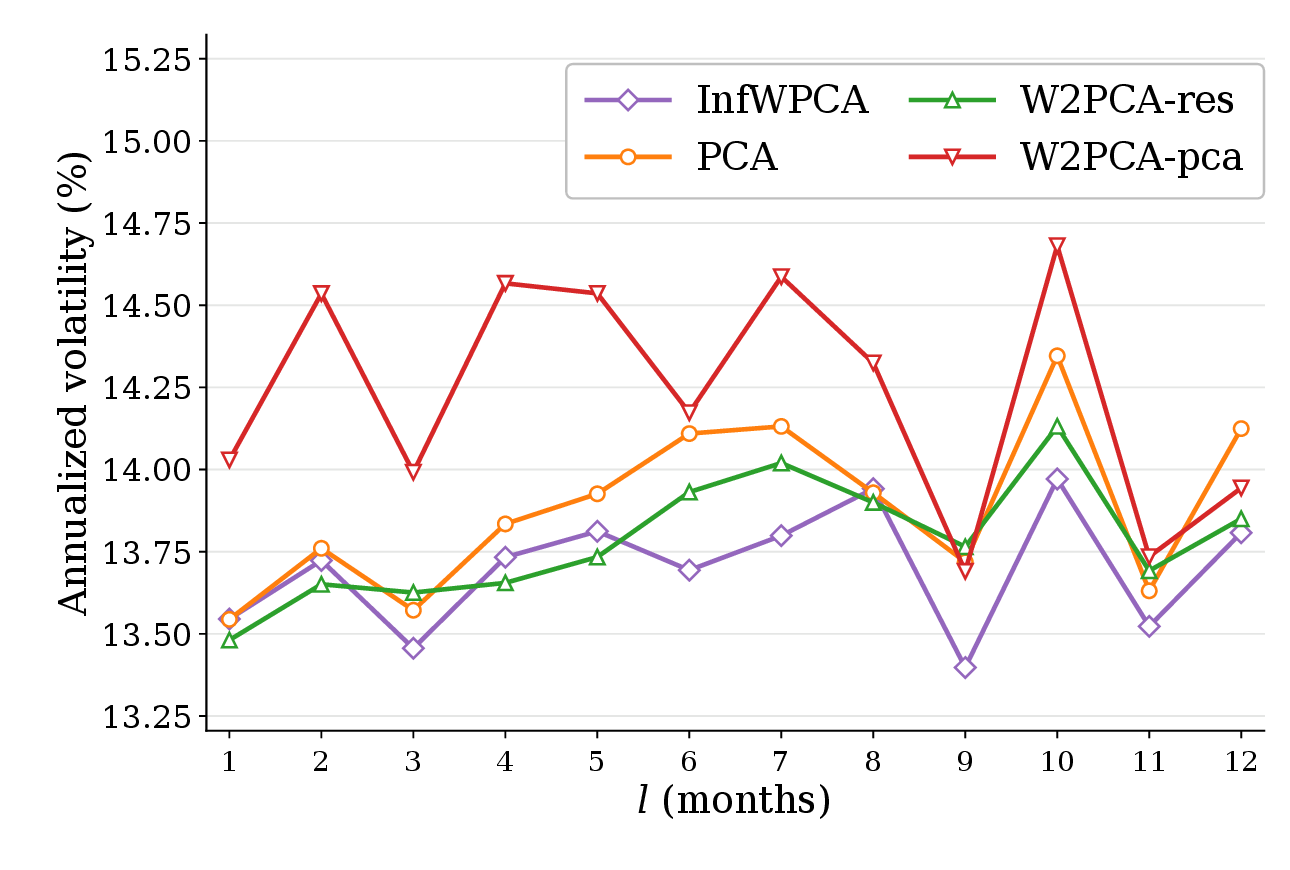}
	}
	\subfloat[Maximum Drawdown]{
		\label{subfig:port_md}
		\includegraphics[width=0.48\linewidth]{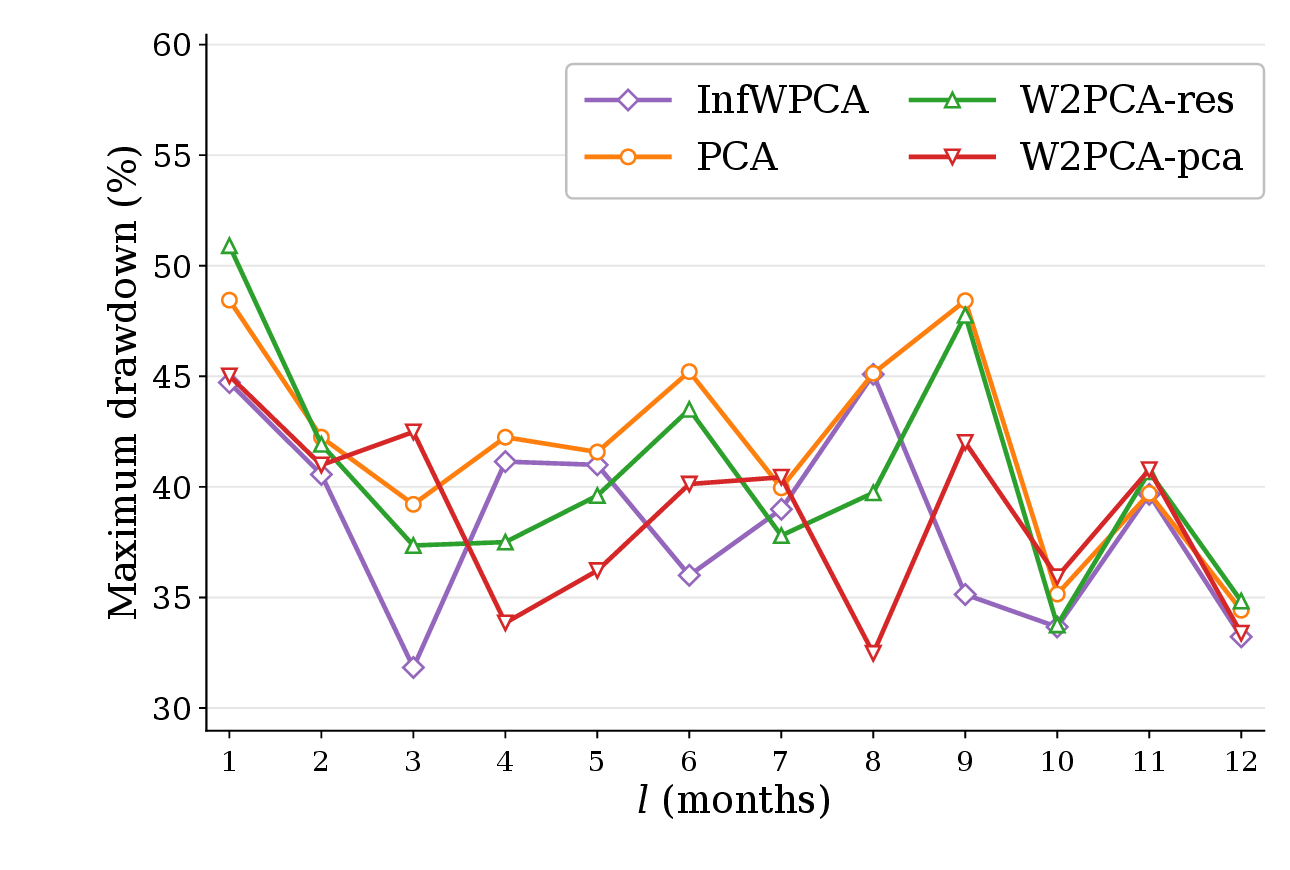}
	}
	
	\caption{Four performance metrics of the tested models under different update intervals.}
	\label{fig:portfolio}
\end{figure}

Table~\ref{tab:portfolio} summarizes the average performance across all twelve update frequencies during the period from January 1974 to December 2025. 
We can observe that InfWPCA outperforms the competing models on all four metrics. 
Figure~\ref{fig:portfolio} provides a more detailed view of its performance under different update intervals. 
InfWPCA attains the highest Sharpe ratio at nine of the twelve update intervals and the highest annualized mean return at seven. 
Moreover, it records the lowest annualized volatility and maximum drawdown at seven update intervals each. 
Taken together, these observations reveal a coherent empirical advantage of InfWPCA, suggesting that a more reliable principal subspace can translate into better downstream portfolio allocation decisions.

\section{Conclusion}

\label{sec:conclusion}

At its core, PCA is concerned with a fundamental choice: which directions of variation should be preserved, and which can be discarded. 
When the underlying distribution is not fully known, this choice depends not only on the observed data but also on how uncertainty is modeled. 
A central insight emerging from this work is that distributional robustness need not serve merely as a safeguard against perturbations. 
It can play an active role in determining the principal subspace itself. 
The proposed infinity-Wasserstein PCA formulation makes this possibility precise. 
By retaining information from individual reconstruction residuals beyond empirical covariance, it allows distributional uncertainty to influence subspace selection. 
In this sense, the geometry of transportation is more than a technical ingredient in defining an ambiguity set. 
It also shapes the low-dimensional representation that PCA seeks to learn. 
This structural perspective underpins the exact reformulation, finite-sample guarantee, and efficient algorithm developed in this paper. 
These findings open several avenues for further investigation. 
One is the possibility of incorporating additional structural requirements into the proposed framework, such as sparsity \cite{Wang2025Decentralized} and nonnegativity \cite{Wang2025Support}. 
Another promising direction is to integrate principal subspace estimation with downstream decision-making \cite{Jiang2026Optimized,Kao2014Directed}. 
Rather than treating the solution of PCA as a fixed input, such a task-aware approach would enable specific applications to directly guide the selection of principal subspace. 

\bibliographystyle{abbrv}
\bibliography{InfWPCA-library}

\addcontentsline{toc}{section}{References}

\end{document}